\documentclass[a4paper,11pt]{article}

\usepackage{amssymb,latexsym}
\usepackage{euscript}
\usepackage{eufrak}
\usepackage{amsmath}
\usepackage{tikz}
\usepackage{bussproofs}

\newtheorem{defi}{\bf Definition}[section]
\newtheorem{rem}{\bf Remark}[section]
\newtheorem{lem}{\bf Lemma}[section]
\newtheorem{prop}{\bf Proposition}[section]
\newtheorem{teo}{\bf Theorem}[section]
\newtheorem{cor}{\bf Corollary}[section]
\newtheorem{eje}{\bf Examples}[section]

\newenvironment{dem}{\noindent\bf Proof. \rm}{\hfill $\mbox{\boldmath{$ \square$}}$}

\makeatother

\makeatletter
\newcommand*{\mirrormodels}{%
  \mathrel{%
    \mathpalette\reflectmathsymbol\models
  }%
}
\newcommand*{\reflectmathsymbol}[2]{%
  \reflectbox{$\m@th#1#2$}%
}

\makeatletter
\providecommand*{\Dashv}{%
  \mathrel{%
    \mathpalette\@Dashv\models
  }%
}
\newcommand*{\@Dashv}[2]{%
  \reflectbox{$\m@th#1#2$}%
}
\makeatother

\title{ {\bf  On the logic that preserves degrees of truth associated to involutive Stone algebras}}

\author{
    Liliana M. Cant\'u $^\textup{\scriptsize a}$
    \and
    Mart\'{\i}n Figallo
  $^\textup{\scriptsize b}$
}

\date{
    $^\textup{\scriptsize a}$\textit{\small Universidad Tecnol\'ogica Nacional, Facultad Regional Tierra del Fuego -- IPES Florentino Ameghino, Ushuaia, Argentina.}
    \\
    $^\textup{\scriptsize b}$\textit{\small Departamento de Matem\'atica, Universidad Nacional del Sur (UNS), Bah\'ia Blanca, Argentina}
}

\begin{document}

\maketitle

\thispagestyle{myheadings}

\begin{abstract}
Involutive Stone algebras (or {\bf S}--algebras) were introduced by R. Cignoli and M. Sagastume in connection to the theory of $n$-valued \L ukasiewicz--Moisil algebras. 

In this work we focus on the logic that preserves degrees of truth associated to involutive Stone algebras, named {\bf \em Six}. This follows a very general pattern that can be considered for any class of truth structure endowed with an ordering relation, and which intends to exploit many--valuedness focusing on the notion of inference that results from preserving lower bounds of truth values, and hence not only preserving the value $1$.

Among other things, we prove that {\bf \em Six} is a many--valued logic (with six truth values) that can be determined by a finite number of matrices (four matrices). Besides, we show that {\bf \em Six} is a paraconsistent logic. Moreover, we prove that it is a genuine LFI (Logic of Formal Inconsistency) with a consistency operator that can be defined in terms of the original set of connectives. Finally, we study the proof theory of {\bf \em Six} providing a Gentzen calculus for it, which is sound and complete with respect to the logic.

\vspace*{4mm}


\vspace*{4mm}

\end{abstract}

\section{\large \bf Introduction}\label{s1}

Let $A$ be a De Morgan algebra. Denote by $K(A)$ the set of all elements $a\in A$ such that the De Morgan negation $\neg a$ of $a$ coincides with the complement of $a$. For every $a\in A$, let $K_a=\{k\in K(A) : a\leq k\}$ and, if $K_a$ has a least element, denote it by $\nabla a$ the least element of $K_a$. The class of all De Morgan algebras $A$ such that for every $a\in A$, $\nabla a$ exists, and the map 

$$a \mapsto \nabla a$$ 

\noindent is a lattice--homomorphism is called the class of involutive Stone algebras, denoted by $\bf  S$. These algebras were introduced by R. Cignoli and M. Sagastume (\cite{CS1}) in connection to the theory of $n$-valued \L ukasiewicz--Moisil algebras.

Most logics studied in the literature in association with a given class of algebras are defined by taking 1 as the only truth value to
be preserved by inference. These truth-preserving logics do not take full advantage of being many-valued, as they focus on
the truth value 1 (the truth) and not on other intermediate truth values. 

A different approach to preserving only truth is going beyond this by considering the notion of inference that results from preserving lower bounds of truth values, and hence not only preserving the value 1. In this setting the language remains the same as in truth-preserving logics, and what changes is the inference relation. This kind of inference corresponds to the so-called logics preserving degrees of truth discussed at length and follows a very general pattern which could be considered for any class of truth structures endowed with an ordering relation.

The study of logics that preserve degrees of truth goes back to W\'ojcicki in his 1988 book \cite{W}, in the context of \L ukasiewicz logic, and then expanded in \cite{F1, F2, F3, F4, BEFGG}, among others. Regarding these logics and their relation to paraconsistency, it is worth citing \cite{CE, CEG, EFGN}.

In this work, we focus on the logic $\mathbb{L}_{\bf S}^{\leq}$ that preserves degrees of truth associated to involutive Stone algebras. Among other things, we prove that $\mathbb{L}_{\bf S}^{\leq}$ is a many--valued logic (with six truth values) that can be determined by a finite number of matrices (four matrices). Besides, we show that $\mathbb{L}_{\bf S}^{\leq}$ is a paraconsistent logic. Moreover, we prove that it is a genuine LFI (Logic of Formal Inconsistency, \cite{WCMCJM}) with a consistency operator that can be defined in terms of the original set of connectives. Finally, we study the proof theory of $\mathbb{L}_{\bf S}^{\leq}$ providing a Gentzen calculus for it, which is sound and complete with respect to the intended semantics.

\section{\large \bf Preliminaries}\label{s2}

In this section, we shall recall some well-known notions and results in order to make this work self-contained.
Recall that a {\em De Morgan algebra} is a structure $\langle A, \wedge, \vee, \neg, 0, 1 \rangle$ of type $(2,2,1,0,0)$ such that the reduct $\langle A, \wedge, \vee,  0, 1 \rangle$ is a bounded distributive lattice and $\neg$ satisfies the identities:
\begin{itemize}
\item[](DM1) $\neg \neg x \approx x$,
\item[](DM2) $\neg(x \wedge y) \approx \neg x \vee \neg y$.
\end{itemize}

Note that it follows that $\neg(x \vee y) = \neg x \wedge \neg y$, $\neg 0 = 1$, and $\neg 1 = 0$. We denote by $\bf M$ the variety of all De Morgan algebras. Let $A$ be a De Morgan algebra. A set  $F\subseteq A$ is a lattice--filter of $A$ iff  (a) $1\in F$, (b) $x\in F$ and $x \leq y$ imply $y\in F$, and (c) $x, y \in F$ implies $x\wedge y \in F$. By ${\cal F}i(A)$ we denote the family of all lattice--filters of $A$. 

Let $A\in {\bf M}$, $a \in A$ is said to be {\em complemented} if there exists $b\in A$, $a \vee b = 1$ and $a \wedge b = 0$, $b$ is the {\em complement} of $a$ and it is denoted $a'$. If the complement of an element $a$ exists this is unique. We denote by $B(A)$ the {\em center of} $A$, that is, the set of all complemented elements of $A$.  It can be seen that, $B(A)$ is a {\bf M}-subalgebra of $A$.

Let $K(A)$ be the subset of $A$ formed by all elements $a \in B(A)$ such that the De Morgan negation $\neg a$ coincides with the complement $a'$. That is, 

$$K(A)=\{a\in B(A): \neg a = a'\}.$$
    
It is easy to check that  $K(A)$ is a {\bf M}-subalgebra de $B(A)$. 
For every $a\in A$, let $K_a =\{k \in K(A): a\leq k \}$. If $K_a$ has a least element, we denote it by $\nabla a$ as in \cite{CS1}. If for every $a\in A$, $K_a$ has a least element, it is said that $A$ is a $\nabla$-De Morgan algebra. The {\em class of involutive Stone algebras}, denoted by {\bf S}, is the subclass of $\nabla$-De Morgan algebras with the property that the map $a \mapsto \nabla a$ from $A$ into $K(A)$ is a lattice-homomorphism.

\begin{teo}\label{TeoCignoli} (Cignoli et al. (\cite{CS2})) $\bf S$ is an equational class. Indeed, a De Morgan algebra  $A \in {\bf S}$ if and only if there is an operator $\nabla : A\to A$  satisfying the following equations:
\begin{itemize}
\item[(IS1)] $\nabla 0 \approx 0$,
\item[(IS2)] $a \wedge \nabla a \approx a$,
\item[(IS3)] $\nabla (a \wedge b) \approx \nabla a \wedge \nabla b$,
\item[(IS4)] $\neg \nabla a \wedge \nabla a \approx 0$.
\end{itemize} 
\end{teo}

If $A\in {\bf S}$, then $a^\star = \neg \nabla a$ is the pseudocomplement of $A$ (see \cite{V}) and $a^{\star\star}=\neg \nabla \neg \nabla a =\nabla a$. Since $a^{\star} \vee a^{\star\star} = 1$, $A$ is a {\em Stone lattice}. The dual pseudocomplement $a^{+}$ of $a$ also exists, in fact, $a^{+}= \nabla \neg a$ and therefore $A$ is also a {\em dual Stone lattice}. From this, we have that {\bf S} is a subclass of the class of {\em double Stone algebras} (see \cite{CS2}). \\[2mm]
It is well known that in a Stone lattice $A$, $a\in B(A)$ iff $a=a^{\star\star}$. Consequently, if $A$ is a De Morgan algebra that is also a Stone lattice and, if we define $\nabla a=a^{\star\star}$, we have that: the algebra $\langle A,\vee,\wedge,\neg,\nabla,0,1\rangle \in {\bf S}$ \, iff \, $a^{\star}=\neg(a^{\star\star})$ for each $a\in A$  iff  $B(A)=K(A)$  (cf.~\cite{CS2}). \\
Motivated by this, the members of {\bf S} were called {\em involutive Stone algebras}. \\[2mm]
Some examples of involutive Stone algebras are the following. 

\begin{eje}
\begin{itemize}
  \item[]
  \item[(i)] Every Boolean algebra $\langle A, \wedge, \vee, \neg, 0,1\rangle$ admits a structure of involutive Stone algebra. Indeed, taking $\nabla x =x$ we have that $\langle A, \wedge, \vee, \neg, \nabla, 0,1\rangle$ is a involutive Stone algebra. 

  \item[(ii)] Let $\langle A, \wedge, \vee , \sim , \sigma^n_1,...,\sigma^n_{n-1},0,1\rangle$ be an $n$-valued \L ukasiewicz-Moisil algebra. Then, the reduct $\langle A, \wedge, \vee , \sim , \sigma^n_{n-1},0,1 \rangle$ is an involutive Stone algebra (\cite{CS1}). In particular, when $n=3$ we have that $\langle A, \wedge, \vee , \sim , \sigma^3_2,0,1 \rangle \in {\bf S}$ and that the operator $\sigma^3_2$ satisfies the relation \, $a \wedge \sigma^3_2 \neg a \leq b \vee \neg \sigma^3_2 b$ \, for each $a,b \in A$. On the other hand, if $A\in {\bf S}$ and the relation $a\wedge \nabla \neg a \leq b\vee \neg \nabla b$ holds in $A$, then $\langle A,\vee, \wedge, \neg,\Delta,\nabla,0,1\rangle$ is a three-valued \L ukasiewicz algebra, where $\Delta =\neg \nabla \neg$. Therefore, (Cignoli et al.) the class of three-valued \L ukasiewicz algebras coincides with the class of involutive Stone algebras that satisfy the equation $a\wedge \nabla \neg a \leq b\vee \neg \nabla b$.

 \item[(iii)] For each $n\geq 2$, ${\bf \L}_{n}$ will denote the $n$-element chain \, $0 < \frac{1}{n-1} < \dots < \frac{n-2}{n-1} < 1$ with the well-known lattice structure and \, $\neg (\frac{i}{n-1}) = \frac{n-i-1}{n-1}$ where $i=0,1,\dots,n-1$.

The algebras ${\bf \L}_{n}$ are  examples of involutive Stone algebras, where $\nabla x= 1$, for all $x\not= 0$ and $\nabla 0=0$.

 \item[(iv)] The following is another important example of an involutive Stone algebra, as we shall see. Let $\mathbb{S}_6=\{0, \frac{1}{3}, N, B, \frac{2}{3}, 1\}$ with the lattice structure given by

\

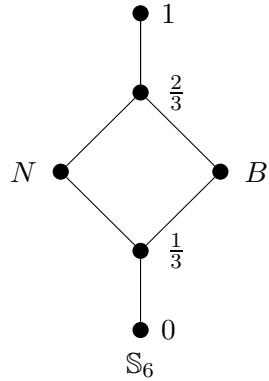
\begin{center}
\begin{tikzpicture}[scale=.7]
\tikzstyle{every node}=[draw,circle,fill=black,inner sep=2pt]

	\node (one) at (4,4.5) [label=right:$1$] {};
	\node (c) at (4,3) [label=right:$\frac{2}{3}$] {};
  \node (b) at (2.5,1.5) [label=left:$N \, $] {};
	 \node (n) at (5.5,1.5) [label=right:$\, B$] {};
 \node (a) at (4,0) [label=right:$\frac{1}{3}$] {};
\node (zero) at (4,-1.5) [label=right:$0$, label=below:$\mathbb{S}_6$] {};

  \draw (zero) -- (a) -- (n) -- (c) -- (one) -- (c) -- (b) -- (a);

\end{tikzpicture} 
\end{center}

the negation given by $\neg 0 =1$, $\neg 1 =0$, $\neg N =N$, $\neg B =B$, $\neg \frac{1}{3} =\frac{2}{3}$ and $\neg \frac{2}{3}=\frac{1}{3}$. Besides,  $\nabla x= 1$, for all $x\not= 0$ and $\nabla 0=0$. 

\end{itemize}
\end{eje}

\begin{rem} The operator $\nabla$ enjoys well--known modal properties such as the Necessitation Rule.
\end{rem}

\

In \cite{CS2}, the authors extended the well--known Priestley's duality for De Morgan algebras (by Cornish and Fowler) to involutive Stone algebras. Using this duality, they determined the subdirectly irreducible and simple {\bf S}--algebras.

\ 

In this work we focus on the logic that preserves degrees of truth associated to involutive Stone algebras. This follows a very general pattern that can be considered for any class of truth structure endowed with an ordering relation; and which intends to exploit many--valuedness focusing on the notion of inference that results from preserving lower bounds of truth values, and hence not only preserving the value $1$ (see \cite{F1,F2,F3}). More precisely, if $K$ is a class of ordered algebras and $\mathfrak{Fm}$ is the absolutely free algebra of the same type of $K$, the notion of inference $\models^{\leq}_{K}$ that results from preserving lower bounds of truth values gives rise to the following definition: for all $\Gamma \cup \{\alpha\}\subseteq Fm$ \, $\Gamma \models^{\leq}_{K} \alpha$ \, iff \,  $\forall A\in {\bf K}$, $\forall v\in Hom_{\bf K}(\mathfrak{Fm},A)$, $\forall a\in A$, \, if $v(\alpha_i)\geq a$, for all $i\leq n$, then $v(\alpha)\geq a$.

\section{\large \bf The Semantic Models: involutive Stone algebras}\label{s3}

In this section, we shall exhibit the structure of our semantic models as well as their most important properties.

\begin{prop}\label{prop1} If $A$ is an involutive Stone algebra then
\begin{itemize}
  \item[(i)] $\nabla A= K(A)$, 
  \item[(ii)] $x\in \nabla A$ \, iff \, $x = \nabla x$,
	\item[(iii)] $\nabla (x \vee y)= \nabla x \vee \nabla y$.
\end{itemize}
\end{prop}
\begin{dem}  It is consequence of Theorem \ref{TeoCignoli}. 
\end{dem}

\

Then, $\nabla$ satisfies the following properties:

\begin{lem}\label{propcuadrado} In every involutive Stone algebra $A$ the following identities hold:

\begin{tabular}{llllll}
(IS5) & $\nabla 1 \approx 1$,\hspace{2cm} & (IS6) & $\neg x \vee \nabla x \approx 1$, & (IS7) & $\nabla \nabla x \approx \nabla x$, \\
(IS8) & $\nabla \neg \nabla x \approx \neg \nabla x$, & (IS9) & $\nabla( x \vee \neg x) \approx 1$,& &\\
(IS10) & $x \wedge \neg \nabla x \approx 0$,& (IS11) & $\nabla(x \vee \nabla y) \approx \nabla x \vee \nabla y$.& &\\
\end{tabular}
\end{lem}
\begin{dem} It is consequence of Theorem \ref{TeoCignoli}. 
\end{dem}

\

As usual, we can define the operator $\Delta$ as follows

$$ \Delta a =_{def.} \neg \nabla \neg a $$

From this definition and Proposition \ref{prop1} (i) we have that \, $\Delta A= \nabla A= K(A)$. Besides, 

\begin{lem} In every involutive Stone algebra $A$ the following identities hold:

\begin{tabular}{llll}
(IS12) &  $\Delta x \wedge  x \approx \Delta x$, & (IS13) & $\Delta \nabla x \approx \nabla x$,\\
(IS14) & $\nabla \Delta x \approx \Delta x$, \hspace{2cm} & (IS15) & $\Delta(x \wedge \neg x)\approx 0$, \\
(IS16) & $\Delta(x \vee y) \approx \Delta x \vee \Delta y$,& (IS17) & $\Delta(x \wedge y) \approx \Delta x \wedge \Delta y$. \\
\end{tabular}
\end{lem}
\begin{dem} It is consequence of Theorem \ref{TeoCignoli}. 
\end{dem}

\

Taking advantage of the topological duality developed in \cite{CS2}, the authors proved the following.

\begin{lem}(Cignoli et al. \cite{CS2}) The subdirectly irreducible algebras of {\bf S} are ${\bf \L}_n$, for $2\leq n\leq5$, and $\mathbb{S}_6$. The simple algebras of {\bf S} are ${\bf \L}_2$ and ${\bf \L}_3$.
\end{lem}

And from the fact that the algebras ${\bf \L}_n$, for $2\leq n\leq5$ are {\bf S}--subalgebras of $\mathbb{S}_6$ we have:

\begin{cor}\label{coro1} $\mathbb{S}_6$ generates the variety {\bf S}.
\end{cor}

\

This last results indicates that a given equation holds in every involutive Stone algebra iff it holds in $\mathbb{S}_6$. This will be very useful in what follows. 

\begin{rem} Belnap's four-valued logic $FOUR$ ([3, 4]), is a logical system that is well known for many applications it has found in several fields. This logic intends to model the following situation: Faced with a situation where one has several conflicting pieces of information on the truth of a sentence, or where
one has no information about it, the classical truth-values (true and false) must be treated as being mutually independent, thus giving birth to four non-classical
epistemic values: 1 (true and not false), 0 (false and not true), {\bf N} (neither true nor false) and  {\bf B} (both true
and false). This gives rise to the well--known four--element De Morgan algebra $\mathfrak{B}_4$, where

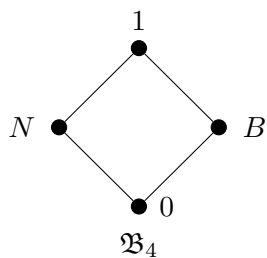
\begin{center}
\begin{tikzpicture}[scale=.7]
\tikzstyle{every node}=[draw,circle,fill=black,inner sep=2pt]
 	
	\node (one) at (-1,3) [label=above:$1$] {};
  \node (b) at (-2.5,1.5) [label=left:$N \, $] {};
	\node (a) at (0.5,1.5) [label=right:$\, B$] {};
  \node (zero) at (-1,0) [label=right:$0$, label=below:$\mathfrak{B}_4$] {};

  \draw (zero) -- (a) -- (one) -- (b) -- (zero);
	\end{tikzpicture}
\end{center}

$\neg 0=1$, $\neg 1 =0$, $\neg N=N$ and $\neg B=B$. \\
On the other hand, many--valued logics (\L ukasiewicz $n$-valued logics, for instance) intend to express the existence of various degrees of truth.
Since ${\bf \L}_4$ and $\mathfrak{B}_4$ are both {\bf M}-subalgebras of $\mathbb{S}_6$, we think that any logic based on this particular algebra would combine this two ideas. In the next section we shall study one of them.  

The lattice $\mathfrak{B}_4$ gives rise to another well-known logic, namely, the {\bf tetravalent modal logic} ${\cal TML}$ which was studied in depth in  \cite{FR} and \cite{MF}.
\end{rem}

\

\section{\large \bf The logic that preserves degrees of truth with respect to {\bf S}}\label{s4}

We now focus on the logic that preserves degrees of truth $\mathbb{L}_{\bf S}^{\leq}$ associated to involutive Stone algebras. 

\

Let $\mathfrak{Fm}=\langle Fm, \wedge, \vee, \neg, \nabla, \bot, \top \rangle$ be the absolutely free algebra of type $(2,2,1,1,0,0)$ generated by some fixed denumerable set of propositional variables $Var$. As usual, the letters $p,q, \dots$ denote propositional variables, the letters $\alpha,\beta, \dots$ denote formulas and $\Gamma,\Sigma, \dots$ sets of formulas. If $A$ and $B$ are to algebras of type $(2,2,1,1,0,0)$ (not necessarily {\bf S}--algebras), we denote by $Hom(A,B)$ the set of all homomorphisms from $A$ into $B$.

\begin{defi} \label{logLogPreserv} The logic that preserves degrees of truth over the variety ${\bf S}$ is $\mathbb{L}^{\leq}_{\bf S}=\langle Fm, \models^{\leq}_{\bf S}\rangle$ where, for every $\Gamma\cup\{\alpha\}\subseteq Fm$: 

\begin{itemize}
\item[(i)] If $\Gamma=\{\alpha_1,\dots, \alpha_n\}\not=\emptyset$,

\begin{center}
\begin{tabular}{ccc}
\hspace{-1cm}$\alpha_1,\dots, \alpha_n \models^{\leq}_{\bf S} \alpha$ & iff &  $\forall A\in {\bf S}$, $\forall v\in Hom(\mathfrak{Fm},A)$, $\forall a\in A$\\
& & if $v(\alpha_i)\geq a$, for all $i\leq n$,\\
& &  then $v(\alpha)\geq a$ \\
\end{tabular}
\end{center}

\item[(ii)] $\emptyset \models^{\leq}_{\bf S} \alpha$ \,  iff \, $\forall A\in {\bf S}$, $\forall v\in Hom(\mathfrak{Fm},A)$, $v(\alpha)=1$.

\item[(iii)] If $\Gamma\subseteq Fm$ is non--finite,

\begin{center}
\begin{tabular}{ccc}
\hspace{-1cm}$\Gamma \models^{\leq}_{\bf S} \alpha$ & iff &  there are $\alpha_1,\dots, \alpha_n \in \Gamma$ such that \\ 
& & $\alpha_1,\dots, \alpha_n \models^{\leq}_{\bf S} \alpha$.\\

\end{tabular}
\end{center}
\end{itemize}
\end{defi}

\

From the above definition we have that $\mathbb{L}^{\leq}_{\bf S}$ is a sentential logic, that is, $\models^{\leq}_{\bf S}$ is a finitary consequence relation over $Fm$.

\begin{prop}\label{lemtripleequiv} Let $\{\alpha_1,\dots, \alpha_n, \alpha\}\subseteq Fm$, $n\geq 1$. Then, the following conditions are equivalent:
\begin{itemize}
\item[(i)] $\alpha_1,\dots, \alpha_n \models^{\leq}_{\bf S} \alpha$,
\item[(ii)] $\alpha_1 \wedge\dots \wedge \alpha_n \models^{\leq}_{\bf S} \alpha $,
\item[(iii)] ${\bf S} \models \alpha_1 \wedge\dots \wedge \alpha_n  \preccurlyeq \alpha$, that is, the inequality  $\alpha_1 \wedge\dots \wedge \alpha_n  \preccurlyeq \alpha$ holds in the variety {\bf S}.
\end{itemize}
\end{prop} 
\begin{dem} (i)$\Rightarrow$(ii) Let $\{\alpha_1,\dots, \alpha_n, \alpha\}\subseteq Fm$, $n\geq 1$.  Suppose that $\alpha_1,\dots, \alpha_n \models^{\leq}_{\bf S} \alpha$. Then, $\forall A\in {\bf S}$, $\forall v\in Hom(\mathfrak{Fm},A)$, $\forall a\in A$, if $v(\alpha_i)\geq a$, for all $i\leq n$, then $v(\alpha)\geq a$.
Since $A$ has a lattice structure, we have that $v(\alpha_i)\geq a$, for all $i\leq n$ iff $\bigwedge_{i=1}^{n}v(\alpha_i)\geq a$ iff $v\left(\bigwedge_{i=1}^{n}(\alpha_i)\right)\geq a$. Then, if $v\left(\bigwedge_{i=1}^{n}(\alpha_i)\right)\geq a$ we have $v(\alpha)\geq a$; that is,  $\alpha_1 \wedge\dots \wedge \alpha_n  \preccurlyeq \alpha$ holds in {\bf S}. Analogously, it is proved (ii)$\Rightarrow$(i). Finally, (iii) is just a different way of writing (ii), assuming (i). 
\end{dem}

\

\begin{rem} The equivalence (i)$\Leftrightarrow$(ii) says that $\mathbb{L}^{\leq}_{\bf S}$ is a {\em conjunctive logic}. On the other hand, (ii)$\Leftrightarrow$(iii) expresses that the consequence relation  $\models^{\leq}_{\bf S}$ internalize the order in {\bf S}.
\end{rem}

\

As a consequence of Proposition \ref{lemtripleequiv} we have the next.

\begin{cor}\label{CorEqSemIgualdad} Let $\alpha,\beta \in Fm$. Then,
\begin{itemize}
\item[(i)] $\alpha \models^{\leq}_{\bf S} \beta$ \, iff \, ${\bf S} \models \alpha \preccurlyeq \beta $,
\item[(ii)] $\alpha \Dashv \models^{\leq}_{\bf S} \beta$ \,  iff  \, ${\bf S} \models \alpha \approx \beta $.
\end{itemize}
\end{cor}
\begin{dem} (i) Immediate taking $n=1$. (ii) Immediate from (i).
\end{dem}

\

From, Corollary \ref{coro1} we have:

\begin{prop} \label{S=S6} For every \ $\Gamma \cup\{\alpha\} \subseteq Fm$,
$\Gamma \models^{\leq}_{\bf S} \alpha$ \ if and only if for every $h
\in Hom(\mathfrak{Fm}, \mathbb{S}_{6})$, $\bigwedge \{
h(\gamma) \ : \ \gamma \in \Gamma\} \leq h(\alpha)$. In particular,
$\emptyset \models^{\leq}_{\bf S} \alpha$ if and only if $h(\alpha)=1$ for all $h \in Hom(\mathfrak{Fm},
\mathbb{S}_{6})$.
\end{prop}

\

\section{$\mathbb{L}^{\leq}_{\bf S}$ as a matrix logic}\label{s5}

In this section we shall see that $\mathbb{L}^{\leq}_{\bf S}$ can be characterized in terms of matrices and generalized matrices (g-matrices).

Recall that a {\em logic matrix} (or simply a matrix) is a pair $\langle A, F\rangle$ where $A$ is an algebra and $F$ is a non--empty subset of $A$. A {\em generalized matrix} is a pair $\langle A, {\cal C}\rangle$ where $A$ is an algebra and ${\cal C}$ is a family of subsets of $A$ which constitutes a closed set system of subsets of $A$. That is, ${\cal C}$ verifies:

\begin{itemize}
\item[(C1)] $A\in {\cal C}$,
\item[(C2)] If $\{X_i\}_{i\in I}\subseteq {\cal C}$ then $\bigcup_{i\in I} X_i \in {\cal C}$.
\end{itemize}

Let $M=\langle A, F\rangle$ be a matrix. The logic  $\mathbb{L}_{M}=\langle \mathfrak{Fm}, \models_{M}\rangle$ determined by $M$ is defined as follows: let $\Gamma\cup\{\alpha\}\subseteq Fm$

\begin{center}
\begin{tabular}{ccc}
$\Gamma \models_{M} \alpha$ & iff &  $\forall v\in Hom(\mathfrak{Fm},A)$\\
& & if $v(\gamma) \in F$, for all $\gamma\in \Gamma$, then $v(\alpha)\in F$. \\
\end{tabular}
\end{center}

If ${\cal F}=\{M_i\}_{i\in I}$ is a family of matrices, then the logic $\mathbb{L}_{\cal F}=\langle \mathfrak{Fm}, \models_{\cal F}\rangle$ determined by the family ${\cal F}$ is the defined by \, $\models_{\cal F} = \bigcap \limits_{i\in I} \models_{M_i}$.

On the other hand, the logic determined by a g-matrix $G=\langle A, {\cal C}\rangle$ is $\mathbb{L}_{G}=\langle \mathfrak{Fm}, \models_{G}\rangle$ where: let $\Gamma\cup\{\alpha\}\subseteq Fm$

\begin{center}
\begin{tabular}{ccc}
$\Gamma \models_{G} \alpha$ & iff &  $\forall v\in Hom(\mathfrak{Fm},A)$, $\forall F\in{\cal C}$\\
& & if $v(\gamma) \in F$, for all $\gamma\in \Gamma$, then $v(\alpha)\in F$. \\
\end{tabular}
\end{center}

The logic determined by a family of g-matrices is defined in a similar way to the logic determined by a family of matrices.

Now, consider the following families of matrices and g-matrices and their associated logics. Let 

$${\cal F}=\{\langle A, F\rangle: A\in {\bf S} \mbox{ and } F \mbox{ is a lattice--filter of } A \}$$

and

$${\cal F}_g=\{\langle A, {\cal F}i(A)\rangle: A\in {\bf S}\}$$

where ${\cal F}i(A)$ is the family of all lattice--filters of $A$. Then,

\begin{teo}\label{teomatrix} The  logic $\mathbb{L}^{\leq}_{\bf S}$ coincides with the logic $\mathbb{L}_{\cal F}$ and also with the logic $\mathbb{L}_{{\cal F}_g}$.
\end{teo}
\begin{dem} The fact that $\mathbb{L}_{\cal F}$ and $\mathbb{L}_{{\cal F}_g}$ coincide is a direct consequence of the definition of logic determined by a family of matrices and g--matrices, respectively. Besides, since ${\bf S}$ is a variety and the notion of filter is elementary definable, the class ${\cal F}_g$ is closed under ultraproducts. Then, $\mathbb{L}_{{\cal F}_g}$ is finitary.

Then, to prove that $\mathbb{L}_{{\cal F}_g}$ and $\mathbb{L}^{\leq}_{\bf S}$ coincide we just must prove it using a finite number of premises.
Let $\alpha_1,\dots,\alpha_n,\alpha\in Fm$ such that $\alpha_1,\dots, \alpha_n \models^{\leq}_{\bf S} \alpha$. Then, by Lemma \ref{lemtripleequiv}, for all $v\in Hom(Fm,A)$, \, (1) $v(\alpha_1) \wedge \dots \wedge(\alpha_n) \leq v(\alpha)$. Let $F\in Fi(A)$ and suppose that $v(\alpha_1), \dots,v(\alpha_n) \in F$, then (by the definition of filter) we have  $v(\alpha_1)\wedge \dots \wedge v(\alpha_n) \in F$. By (1) and $F$ filter, $v(\alpha)\in F$. Then, $\alpha_1,\dots, \alpha_n \models_{{\cal F}_g} \alpha$.

Conversely, suppose that (2) $\alpha_1,\dots, \alpha_n \models_{{\cal F}_g} \alpha$ and let $v\in Hom(Fm,A)$. Consider the filter generated by $v(\alpha_1), \dots, v(\alpha_n)$, $F=Fi(v(\alpha_1), \dots, v(\alpha_n))$. Clearly  $v(\alpha_i)\in F$ for all $i\leq n$, and then by (2), $v(\alpha)\in F$. Then, by properties of generated filters, $v(\alpha_1)\wedge \dots \wedge v(\alpha_n)\leq v(\alpha)$. Therefore, by Lemma \ref{lemtripleequiv}, $\alpha_1,\dots, \alpha_n \models^{\leq}_{\bf S} \alpha$.
\end{dem}

\

On the other hand, the filters of $\mathbb{S}_6$ are precisely the principal filters generated by each element of $\mathbb{S}_6$.
Then, ${\cal F}i(\mathbb{S}_6)=\{[0), [\frac{1}{3}), [N), [B), [\frac{2}{3}), [1)\}$. Now, consider the following family of matrices 

$${\cal F}^{6}=\{\langle \mathbb{S}_6, [0)\rangle, \langle \mathbb{S}_6, \left[\frac{1}{3}\right)\rangle, \langle \mathbb{S}_6, [N)\rangle, \langle \mathbb{S}_6, [B)\rangle, \langle \mathbb{S}_6, \left[\frac{2}{3}\right)\rangle, \langle \mathbb{S}_6, [1)\rangle \}$$

and the generalized matrix

$${\cal F}^{6}_g=\langle\mathbb{S}_6, {\cal F}i(\mathbb{S}_6)\rangle$$

\begin{teo} The logics $\mathbb{L}_{{\cal F}^{6}}$,  $\mathbb{L}_{{\cal F}^{6}_g}$ and $\mathbb{L}^{\leq}_{\bf S}$ coincide.
\end{teo}
\begin{dem}
Since the logic determined by a finite set of matrices is finitary, again we shall prove that they coincide using a finite number of premises. Besides, it is clear that $\mathbb{L}_{{\cal F}^{6}}$ and $\mathbb{L}_{{\cal F}^{6}_g}$ coincide.
On the other hand, let $\alpha_1,\dots,\alpha_n,\alpha\in Fm$ such that $\alpha_1,\dots, \alpha_n \models^{\leq}_{\bf S} \alpha$. By Proposition \ref{S=S6}, this is equivalent to $\mathbb{S}_6 \models \alpha_1 \wedge\dots \wedge \alpha_n  \preccurlyeq \alpha$, that is, for all $h \in Hom(\mathfrak{Fm}, \mathbb{S}_{6})$, $h(\alpha_1)\wedge \dots \wedge h(\alpha_n)\leq h(\alpha)$. Then, by the same properties of filters used before, we have $\alpha_1,\dots, \alpha_n \models_{{\cal F}^{6}_g} \alpha$. The converse is analogous to the converse showed in the proof of Theorem \ref{teomatrix}.
\end{dem}

\

So far, we have proved that $\mathbb{L}^{\leq}_{\bf S}$ is a matrix logic determined by a family of six finite matrices. Next, we shall see that two of them are superfluous.

In first place, observe that the matrix $\langle \mathbb{S}_6, [0)\rangle= \langle \mathbb{S}_6, \{0,\frac{1}{3},N,B,\frac{2}{3}, 1\}\rangle$ is trivial. Indeed,  for all $\alpha\in Fm$ and all $h \in Hom(\mathfrak{Fm}, \mathbb{S}_{6})$ we have that $h(\alpha)\in [0)$. Then, for every $\Gamma\subseteq Fm$ it holds  $\Gamma \models_{\langle \mathbb{S}_6, [0)\rangle} \alpha$ and so we can discard it.

On the other hand, we shall see that $\langle \mathbb{S}_6, [N)\rangle$ and $\langle \mathbb{S}_6, [B)\rangle$ are isomorphic matrices. For this, the next technical result will be useful.

\

\begin{lem}\label{Lemhh'} Let $h: \mathfrak{Fm} \to \mathbb{S}_{6}$,
${\cal V}' \subseteq Var$ \, and \, $h':\mathfrak{Fm} \to \mathbb{S}_{6}$
such that, for all $p\in{\cal V}'$,\\[4mm]
$h'(p)= \left \{ \begin{tabular}{rl}
$h(p)$ & if $h(p)\in \{0,\frac{1}{3},\frac{2}{3},1\}$\\
$N$ & if $h(p)=B$\\
$B$ & if $h(p)=N$\\
\end{tabular}\right.$ . Then\,  $h'(\alpha)= \left \{ \begin{tabular}{rl}
$h(\alpha)$ & if $h(\alpha)\in \{0,\frac{1}{3},\frac{2}{3},1\}$\\
$N$ & if $h(\alpha)=B$\\
$B$ & if $h(\alpha)=N$\\
\end{tabular}\right.$\\[3mm]
for all $\alpha \in \mathfrak{Fm}$ such that $Var(\alpha)\subseteq {\cal V}'$.
\end{lem}
\begin{dem} We use induction on the complexity of the formula $\alpha$, $c(\alpha)$. First, observe that $\{0,\frac{1}{3},\frac{2}{3},1\}\simeq {\bf \L}_4$ is an {\bf S}-subalgebra of $\mathbb{S}_6$, that is, it is closed under all the operations of {\bf S}-algebras.

If $\alpha$ is $p$, where $p$ is a propositional variable, the lemma holds by hypothesis.

Suppose that the lemma holds for every formula with complexity less than $k$, $k\geq1$; and let $\alpha$ a formula with complexity $k$. Then, we have the following cases:

{\bf Case I :} $\alpha$ is $\neg \beta$, with $\beta\in Fm$. If $h(\alpha)\in \{0,\frac{1}{3},\frac{2}{3},1\}$, since $h(\alpha)=h(\neg \beta)=\neg h(\beta)$, we have that $h(\beta)\in \{0,\frac{1}{3},\frac{2}{3},1\}$. On the other hand, $h'(\beta)\in \{0,\frac{1}{3},\frac{2}{3},1\}$. By I.H, $h(\beta)=h'(\beta)$. Since $h'$ is a homomorphism $h'(\alpha)=h'(\neg\beta)=\neg h(\beta)=\neg h'(\beta)=h'(\neg \beta)=h'(\alpha)$. If $h(\alpha)= N$, from $\neg N=N$, we have that $h(\beta)=N$. By I.H., $h'(\beta)=B$ and since $\neg B=B$, we have that $h'(\alpha)=B$. Analogously, we prove the lemma in case that $h(\alpha)=B$.

{\bf Case II :} $\alpha$ is $\beta_1 \wedge \beta_2$, where $\beta_1,\beta_2\in Fm$. If $h(\alpha)=h(\beta_1) \wedge h(\beta_2)\in \{0,\frac{1}{3},\frac{2}{3},1\}$ there are many sub-cases. If $h(\beta_1), h(\beta_2)\in \{0,\frac{1}{3},\frac{2}{3},1\}$, by I.H., $h(\beta_1)=h'(\beta_1)$ and $h(\beta_2)=h'(\beta_2)$ and then $h(\alpha)=h(\beta_1) \wedge h(\beta_2)= h'(\beta_1) \wedge h'(\beta_2)=h'(\alpha)$. If $h(\beta_1)=N$ and $h(\beta_2)\in \{0,\frac{1}{3}\}$, by I.H., $h'(\beta_1)=B$ and $h'(\beta_2)=h(\beta_2)$ and then $h'(\alpha)= h'(\beta_1) \wedge h'(\beta_2)= B\wedge h(\beta_2)=h(\beta_2)= N \wedge h(\beta_2)=h(\beta_1) \wedge h(\beta_2)=h(\alpha)$. The other sub-cases are studied similarly.

{\bf Case III :} $\alpha$ is $\beta_1 \vee \beta_2$, where $\beta_1,\beta_2\in Fm$. It is analogous to case II.

{\bf Caso IV :} $\alpha$ is $\nabla \beta$, with $\beta\in Fm$. Then, $h(\alpha)=h(\nabla \beta)=\nabla h(\beta)\in\{0,1\}$. If, $h(\alpha)=0$, then $h(\beta)=0$, by the definition of $\nabla$ in $\mathbb{S}_{6}$. By I.H., $h'(\beta)=h(\beta)=0$ and therefore $h'(\alpha)=0=h(\alpha)$. If $h(\alpha)=1$, then $h(\beta)\in  \{\frac{1}{3},N,B,\frac{2}{3}, 1\}$. In any of this cases we have that $h'(\alpha)=1$.
\end{dem}

\

Now, we can see that the matrices  $\langle \mathbb{S}_6, [N)\rangle$ and $\langle \mathbb{S}_6, [B)\rangle$ determine the same logic.

\

\begin{lem} The logics $\models_{\langle \mathbb{S}_6, [N)\rangle}$ and $\models_{\langle \mathbb{S}_6, [B)\rangle}$ coincide.
\end{lem}
\begin{dem} Since $\models_{\langle \mathbb{S}_6, [N)\rangle}$ (and $\models_{\langle \mathbb{S}_6, [B)\rangle}$) is a finitary consequence relation and conjunctive; we only need to consider inferences of the form $\alpha \models_{\langle \mathbb{S}_6, [N)\rangle} \beta$ (if  $\beta$ is logically valid we consider $\alpha$ as $\neg\bot$).
Suppose that (1)$\alpha \models_{\langle \mathbb{S}_6, [N)\rangle} \beta$ and let $h \in Hom(\mathfrak{Fm}, \mathbb{S}_{6})$ such that $h(\alpha) \in
\{B,\frac{2}{3}, 1\}$.
If $h(\alpha)=B$, consider the homomorphism $h'$ as in Lemma \ref{Lemhh'}. Then, $h'(\alpha)=N \in \{N,\frac{2}{3}, 1\}$ and by (1), $h'(\beta)\in \{N,\frac{2}{3}, 1\}$. If $h'(\beta)=N$, then $h(\beta)=B \in \{B,\frac{2}{3}, 1\}$. Otherwise, if $h'(\beta)\in\{\frac{2}{3}, 1\}$ then $h'(\beta)=h(\beta) \in\{\frac{2}{3}, 1\} \subseteq \{B,\frac{2}{3}, 1\}$.

If $h(\alpha)\in \{\frac{2}{3}, 1\}$, by (1), $h(\beta)\in \{N,\frac{2}{3}, 1\}$. If $h(\beta)=N$ then $h'(\beta)=B$. But, in this case, $h(\alpha)=h'(\alpha)\in  \{\frac{2}{3}, 1\}$, that is, $h'\in Hom(\mathfrak{Fm}, \mathbb{S}_{6})$ verifies $h'(\alpha)\in \{N,\frac{2}{3}, 1\}$ but $h'(\beta)\notin \{N,\frac{2}{3}, 1\}$, which contradicts (1). Then, $h(\beta)\in \{\frac{2}{3}, 1\}\subseteq \{B,\frac{2}{3}, 1\}$. 
From analyzing all these cases we conclude that $\alpha \models_{\langle \mathbb{S}_6, [B)\rangle} \beta$.
Analogously, we prove that if $\alpha \models_{\langle \mathbb{S}_6, [B)\rangle} \beta$ then $\alpha \models_{\langle \mathbb{S}_6, [N)\rangle} \beta$.
\end{dem}

\

From the results above mentioned, we have:

\

\begin{teo}\label{TeoSix} The logic that preserves degrees of truth associated to involutive Stone algebras, $\mathbb{L}^{\leq}_{\bf S}$, is a six--valued logic determined by the matrices $\langle \mathbb{S}_6, \left[\frac{1}{3}\right)\rangle$, $\langle \mathbb{S}_6, [N)\rangle$, $\langle \mathbb{S}_6, \left[\frac{2}{3}\right)\rangle$ and $\langle \mathbb{S}_6, [1)\rangle$.
\end{teo}

\

In what follows, we shall denote by {\bf \em Six} the logic $\mathbb{L}^{\leq}_{\bf S}$. 

\

Recall that we say that a set of connectives in a given logic is {\em functionally complete} if every possible truth--table on the matrix semantic, which is sound and complete with respect to the logic, can be expressed  in terms of these connectives. We say that the logic is functionally complete if it has a functionally complete set of connectives.

\begin{lem}\label{Lemhh'2} Let $h: \mathfrak{Fm} \to \mathbb{S}_{6}$, $h':\mathfrak{Fm} \to \mathbb{S}_{6}$ \, and \,  $p \in Var$ such that $h(p)=B$ and $h'(p)= N$. Then, $h(\alpha)\not=N$ and $h'(\alpha)\not=B$ for all  $\alpha \in \mathfrak{Fm}$ such that $Var(\alpha) \subseteq \{p\}$.
\end{lem}
\begin{dem} It is a direct consequence of the fact that $\{0,\frac{1}{3},B, \frac{2}{3},1\}$ and $\{0,\frac{1}{3},N, \frac{2}{3},1\}$ are {\bf S}-subalgebras of $\mathbb{S}_{6}$.
\end{dem}

\

\begin{cor} \label{noFuncComple} {\bf \em Six} is not functionally complete.
\end{cor}
\begin{dem}
By  Lemma~\ref{Lemhh'2}, it is not possible to define the function $f:\mathbb{S}_{6} \to \mathbb{S}_{6}$ such that $f(x)=N$, for all $x$, in terms of the connectives.
\end{dem}

\section{{\bf \em Six} as a paraconsistent logic}

In this section, we shall study {\bf \em Six} under the perspective of paraconsistency. We shall see that contradictions (with respect to $\neg$) not necessarily trivialize the inferences, and therefore, {\bf \em Six} is a paraconsistent logic in the sense of da Costa (cf.~\cite{daC,daC2}). Moreover, we shall see that {\bf \em Six} is a Logic of Formal Inconsistency (LFI); and that it is gently explosive with respect to a proper set of formulas $\bigcirc(p)$ which depends on the propositional variable $p$ (cf.~\cite{Tax,WCMCJM}).

\

\begin{prop} {\bf \em Six} is non--trivial and non--explosive.
\end{prop}
\begin{dem} Let  $p$ and $q$ be two different propositional variables and consider the homomorphism $h: \mathfrak{Fm} \to \mathbb{S}_{6}$ such that $h(p)=N$ and $h(q)=B$. Then, \, $p, \neg p \not \models_{\mbox{\small \bf \em Six}} q$.
\end{dem}

\

Considering the same homomorphism as in the proof above we conclude also that \, $\not \models_{\mbox{\small \bf \em Six}} q \vee \neg q$, and therefore

\begin{prop} {\bf \em Six} is paracomplete.
\end{prop}

\

On the other side, in the language of  {\bf \em Six},  put $\bigcirc(p)=\{ \Delta p \vee \Delta \neg p \}$. Then,

\

\begin{lem}  {\bf \em Six} is finitely gently explosive with respect to $\bigcirc(p)$ and $\neg$.
\end{lem}
\begin{dem} Let $p$ and $q$ be two different propositional variables. It is easy to check that
$$\bigcirc(p), p \not \models_{\mbox{\small \bf \em Six}} q \ \ \textrm{ and } \ \
\bigcirc(p), \neg p \not \models_{\mbox{\small \bf \em Six}} q.$$ But, since $h((\Delta p \vee \Delta \neg p) \wedge
p \wedge \neg p) = h((\neg \nabla \neg p \vee \neg \nabla \neg \neg p) \wedge p \wedge \neg p) = h((\neg \nabla \neg p \vee \neg \nabla p) \wedge p \wedge \neg p)=h((\neg \nabla \neg p  \wedge p \wedge \neg p) \vee (\neg \nabla p\wedge p \wedge \neg p))=0 $ (by (IS10)), for all $h \in Hom(\mathfrak{Fm},
\mathbb{S}_{6})$. So,  \, $\bigcirc(p), p, \neg p \models_{\mbox{\small \bf \em Six}} \bot$.
\end{dem}

\

From all the above we have:

\begin{teo}  {\bf \em Six} is an LFI with respect to $\neg$ and with consistency operator $\circ$ defined by $\circ \alpha = \Delta \alpha \vee \Delta \neg \alpha$, for all $\alpha \in Fm$.
\end{teo}

\

Besides, as it happens in the systems $C_n$ of da Costa (cf.~\cite{daC,daC2}), the operator of consistency propagates through the others connectives.

\begin{teo}\label{propag} In {\bf \em Six} it holds:\\[2mm]
\begin{tabular}{ll}
{\rm(i)} $\models_{\mbox{\small \bf \em Six}} \circ \bot$ &
{\rm(ii)} $\circ \alpha \models_{\mbox{\small \bf \em Six}} \circ \nabla \alpha$ \\[2mm]
{\rm(iii)} $\circ \alpha \models_{\mbox{\small \bf \em Six}} \circ \neg \alpha$ &
{\rm(iv)}  $\circ \alpha, \circ \beta \models_{\mbox{\small \bf \em Six}} \circ (\alpha \# \beta)$ \, for \, $\# \in \{\wedge, \vee\}$.
\end{tabular}
\end{teo}
\begin{dem} It is routine.
\end{dem}

\

Moreover, \ $\models_{\mbox{\small \bf \em Six}} \circ \neg^{n} {\circ}
\alpha$, for all $n\geq 0$. In particular, $\models_{\mbox{\small \bf \em Six}}
{\circ} {\circ} \alpha$. So, {\bf \em Six} validates all axioms $(cc)_{n}$ of the logic {\bf mCi} (see \cite{WCMCJM}). As it is usual in the framework of LFI's, it is also possible to define an inconsistency operator~$\bullet$ \, in {\bf \em Six} \, in the following way:

$$ \bullet \alpha =_{def} \neg \circ \alpha.$$

\

Then, $\bullet \alpha$ is logically equivalent to $\nabla \alpha \wedge \nabla \neg \alpha$. Note that $\bullet$
and $\circ$ can be expressed equivalently as  \, $\bullet \alpha = \nabla( \alpha \wedge \neg \alpha)$ \, and \, $\circ \alpha
= \Delta(\alpha \vee \neg \alpha)$. The truth--table of both connectives is the following:

\

\begin{center}
\begin{tabular}{|c|c|c|}         \hline
$p$ & $\circ p$ & $\bullet p$\\ \hline
0       & 1 & 0         \\ \hline
$\frac{1}{3}$ & 0 & 1 \\ \hline
N       & 0 & 1          \\ \hline
B       & 0 & 1         \\ \hline
$\frac{2}{3}$ & 0 & 1 \\ \hline
1       & 1 & 0         \\ \hline
\end{tabular}
\end{center}

\

\begin{teo} In {\bf \em Six} it holds:

\begin{itemize}
\item[{\rm(i)}] $ \alpha \wedge \neg \alpha \models_{\mbox{\small \bf \em Six}} \bullet \alpha$ \, but \, $\bullet \alpha \not \models_{\mbox{\small \bf \em Six}}  \alpha \wedge \neg \alpha$,
\item[{\rm(ii)}] $\bullet \alpha \models_{\mbox{\small \bf \em Six}} \bullet \neg \alpha$ \, and \, $ \bullet \neg \alpha \models_{\mbox{\small \bf \em Six}}\bullet \alpha$,
\item[{\rm(iii)}] $\bullet (\alpha \# \beta) \models_{\mbox{\small \bf \em Six}} \bullet \alpha \vee \bullet \beta$ \, for \, $\# \in \{\wedge, \vee\}$; the converse does not hold.
\end{itemize}
\end{teo}
\begin{dem} It is routine.
\end{dem}

\

This last result show us that the concept of consistency on one side, and contradiction on the other, can be differentiated in  {\bf \em Six}. This is a very valuable characteristic in the universe of LFI's only satisfied by a few of them such as  {\bf
mbC} and $\cal TML$ (see \cite{MF}), the first being the weakest in the hierarchy presented in~\cite{WCMCJM}.
In spite of the fact that {\bf \em Six}  is not functionally complete (cf. Corollary~\ref{noFuncComple}), it enjoys of a nice expressive power. For instance, in {\bf \em Six} it is possible distinguish the ``classical'' truth values  ($0$ and $1$) from the ``non--classical'' ones ($\frac{1}{3}$, $N$, $B$ and $\frac{2}{3}$).

\

To end this section, we shall show that, in {\bf \em Six}, it is possible to reproduce the classical logic. That is, we shall prove a {\em Derivability Adjustment Theorem} (DAT) with respect to the classical propotitional logic {\bf CPL}.

We can think of {\bf CPL} in the language generated by  $\wedge$,  $\vee$ and $\neg$. Let $\mathfrak{Fm}_{\bf CPL}$ the algebra of formulas of {\bf CPL} in this language.

\

\begin{teo} \label{DAT} Let $\Gamma\cup\{\alpha\}$ be a finite set of formulas in {\bf CPL}.
Then,\\
$\Gamma \vdash_{\rm \bf CPL} \alpha$ \, iff \, $\Gamma, \circ p_1, \dots, \circ p_n \models_{\mbox{\small \bf \em Six}} \alpha$
\,  where ${\rm Var}(\Gamma \cup \{ \alpha\})=\{ p_1, \dots,
p_n\}$.
\end{teo}
\begin{dem}($\Leftarrow$) Let $\Gamma = \{\alpha_1, \dots, \alpha_k\} \subseteq \mathfrak{Fm}_{\bf CPL}$ and let
$h \in Hom(\mathfrak{Fm}_{\bf CPL}, \mathbb{B}_{1})$ where
$\mathbb{B}_{1}$ is the two--element Boolean algebra (that is,
$\mathbb{B}_{1}=\{0,1\}$). Since $\mathbb{B}_{1}$ can be seen as an {\bf S}-subalgebra of $\mathbb{S}_{6}$ (if we put
$\nabla(x)=x$, for all $x$, then $\mathbb{B}_{1}\simeq {\bf \L}_{2}$), we have that que $h \in
Hom(\mathfrak{Fm}, \mathbb{S}_{6})$. Then, by hypothesis we have 

$$h(\bigwedge \limits_{i=1}^{k} \alpha_i \wedge \bigwedge \limits_{j=1}^{n} \circ p_j) \leq h(\alpha).$$

Since $h \in Hom(\mathfrak{Fm}, \mathbb{B}_{1})$, then $h(p_j) \in \{0,1\}$ for all $j$, $1 \leq j
\leq n$. In this way, by definition of $\circ$, $h(\circ p_j)= \circ h(p_j) = 1$ for all $j$, $1 \leq j \leq n$. Then, if
$h(\alpha_i)=1$ for all $i$, $1\leq i \leq k$ we have $h(\bigwedge \limits_{i=1}^{k} \alpha_i \wedge \bigwedge
\limits_{j=1}^{n} \circ p_j)=1$ and therefore  $h(\alpha)=1$.
That is, $\Gamma \vdash_{\rm \bf CPL} \alpha$.

\

\noindent ($\Rightarrow$) Suppose that $\Gamma \vdash_{\rm \bf CPL} \alpha$ and let $h \in Hom(\mathfrak{Fm}, \mathbb{S}_{6})$. Let's see that $\Gamma, \circ p_1, \dots, \circ p_n \models_{M} \alpha$ for each of the matrices $M$ which determine {\bf \em Six}. Let $M=\langle \mathbb{S}_{6},[\frac{1}{3})\rangle$ and suppose that  $h(\Gamma \cup \{\circ p_1, \dots, \circ p_n\}) \subseteq[\frac{1}{3})$, where ${\rm Var}(\Gamma \cup \{ \alpha\})=\{ p_1, \dots,p_n\}$. Then, $h(\circ p_i)\not=0$ which implies that $h(p_i)\in \{1,0\}$ for all $i$. Then, we can think that $h \in Hom(\mathfrak{Fm}_{\bf CPL}, \mathbb{B}_{1})$. 

Since $\Gamma \vdash_{\rm \bf CPL} \alpha$ and $h(\Gamma) \subseteq\{1\}$, then $h(\alpha)=1$ and $h(\alpha) \in [\frac{1}{3})$. Then, $\Gamma, \circ p_1, \dots, \circ p_n \models_{\langle \mathbb{S}_{6},[\frac{1}{3})\rangle} \alpha$.

Analogously it is proved  that  $\Gamma, \circ p_1, \dots, \circ p_n \models_{M} \alpha$ for the others matrices which define {\bf \em Six}. Therefore, $\Gamma, \circ p_1, \dots, \circ p_n \models_{\mbox{\small \bf \em Six}} \alpha$.
\end{dem}

\

\section{Proof theory for {\bf \em Six}}\label{s6}

So far, we have studied the logic {\bf \em Six} by semantical tools. In this section, we shall present a syntactic version of it by means of a sequent calculus. If $\Sigma$ is a set of formulas, then $\nabla \Sigma=\{\nabla \alpha: \alpha\in \Sigma\}$.

\begin{defi} We denote  by  $\mathfrak{S}$  the  sequent calculus whose  axioms  and  rules 
are  the following:

\noindent {\bf Axioms}
$$ \mbox{(structural axiom) \, } \displaystyle {\alpha \Rightarrow \alpha} \hspace{2cm} \mbox{($\bot$) \,  } {\bot\Rightarrow } \hspace{2cm} \mbox{($\top$) \,  } {\Rightarrow \top }$$

$$\mbox{(First modal axiom) \,  } {\alpha\Rightarrow \nabla\alpha} \hspace{2cm} \mbox{(Second modal axiom) \,  } {\Rightarrow \nabla\alpha \vee \neg \nabla \alpha}$$

\noindent {\bf Structural rules}

$$ \mbox{(Left weakening) \, } \displaystyle \frac{\Gamma \Rightarrow \Sigma} {\Gamma, \alpha \Rightarrow \Sigma} \hspace{2cm} \mbox{(Right weakening) \, }
 \displaystyle \frac{\Gamma \Rightarrow \Sigma} {\Gamma \Rightarrow \Sigma, \alpha} $$

$$ \mbox{(Cut) \, }
 \displaystyle \frac{\Gamma \Rightarrow \Sigma, \alpha  \hspace{0.5cm} \alpha, \Gamma \Rightarrow \Sigma}{\Gamma \Rightarrow \Sigma} $$

\noindent {\bf Logic Rules}

$$ \mbox{($\wedge \Rightarrow$) \, } \displaystyle \frac{\Gamma, \alpha, \beta \Rightarrow \Sigma} {\Gamma, \alpha \wedge \beta \Rightarrow \Sigma} \hspace{2cm} \mbox{($\Rightarrow \wedge$) \, }
 \displaystyle \frac{\Gamma \Rightarrow \Sigma, \alpha  \hspace{0.5cm} \Gamma \Rightarrow \Sigma, \beta}{\Gamma \Rightarrow \Sigma, \alpha \wedge \beta} $$

$$ \mbox{($\vee \Rightarrow$) \, } \displaystyle \frac{\Gamma, \alpha \Rightarrow \Sigma \hspace{0.5cm} \Gamma, \beta \Rightarrow \Sigma} {\Gamma, \alpha \vee \beta \Rightarrow \Sigma} \hspace{2cm}  \mbox{($\Rightarrow \vee$) \, }
 \displaystyle \frac{\Gamma \Rightarrow \Sigma, \alpha, \beta }{\Gamma \Rightarrow \Sigma, \alpha \vee \beta} $$

$$ \mbox{($\neg$) \, } \displaystyle \frac{\alpha \Rightarrow \beta } {\neg \beta \Rightarrow \neg \alpha } $$

$$ \mbox{($\neg \neg \Rightarrow$) \, } \displaystyle \frac{\Gamma, \alpha \Rightarrow \Sigma } {\Gamma, \neg \neg \alpha \Rightarrow \Sigma} \hspace{2cm} \mbox{($\Rightarrow \neg \neg$)} \, \frac{\Gamma \Rightarrow \alpha, \Sigma } {\Gamma \Rightarrow \neg \neg \alpha, \Sigma} $$

$$ \mbox{($\nabla$) \, } \displaystyle \frac{\Gamma, \alpha \Rightarrow \nabla \Sigma } {\Gamma, \nabla \alpha \Rightarrow \nabla \Sigma} \hspace{2cm} \mbox{($\neg \nabla \Rightarrow$) \, } \displaystyle \frac{\Gamma, \neg \nabla \alpha \Rightarrow \Sigma } {\Gamma, \nabla \neg \nabla \alpha \Rightarrow \Sigma}   $$

\end{defi}

\

The notion of {\em derivation} (o  $\mathfrak{S}$--{\em proof}) in the Gentzen calculus $\mathfrak{S}$ is the usual. That is, we say  that a sequent  $\Gamma \Rightarrow \Sigma$ is {\em derivable}  in $\mathfrak{S}$, denoted $\mathfrak{S} \vdash \Gamma \Rightarrow \Sigma$,  when it has a derivation  whose initial sequents  are instances  of  the  axioms  and  where the  rules used  are among  those in the above  list. We  have  not explicitly  included  the  rules  of Exchange  and  Contraction  because  we are  using sets of formulas  in our sequents,  and not  just multisets or sequences;  thus this Gentzen calculus satisfies  all structural  rules.
If $\alpha$  is a formula we write $\mathfrak{S}\vdash \alpha$ to indicate that $\mathfrak{S} \vdash \Rightarrow \alpha$.

\begin{lem}\label{Lem4Equiv} Let $\alpha_1,\dots, \alpha_n, \beta_1, \dots, \beta_m \in Fm$. The following conditions are equivalent.
\begin{itemize}
  \item[(i)] $\alpha_1,\dots, \alpha_n \Rightarrow \beta_1, \dots, \beta_m$ is derivable in $\mathfrak{S}$, 
  \item[(ii)] $\alpha_1\wedge \dots \wedge \alpha_n \Rightarrow \beta_1, \dots, \beta_m$ is derivable in $\mathfrak{S}$
  \item[(iii)] $\alpha_1,\dots, \alpha_n \Rightarrow \beta_1\vee \dots \vee \beta_m$ is derivable in $\mathfrak{S}$
  \item[(iv)] $\alpha_1\wedge \dots \wedge \alpha_n \Rightarrow \beta_1\vee \dots \vee \beta_m$ is derivable in $\mathfrak{S}$.
\end{itemize}
\end{lem}
\begin{dem} It is routine taking into account that the rules ($\wedge \Rightarrow$), ($\Rightarrow \wedge$), ($\vee \Rightarrow$), ($\Rightarrow \vee$) are the same rules of the classical case.
\end{dem}

\

\begin{rem}\label{Obs1} Since the logical rules which govern the behaviour of the conjunction $\wedge$ and the disjunction $\vee$ are the classical ones, we know that in $\mathfrak{S}$ all classical properties concerning this two connectives hold. In particular, the following hold: let $\alpha, \beta, \gamma \in Fm$. 
If $\mathfrak{S} \vdash \alpha \Leftrightarrow \beta$ then, 
\begin{itemize}
\item[(i)] If $\mathfrak{S} \vdash \alpha \Leftrightarrow \beta$ then, $\mathfrak{S} \vdash \alpha \vee \gamma  \Leftrightarrow \beta \vee \gamma$ \, and \, $\mathfrak{S} \vdash \alpha \wedge \gamma \Leftrightarrow \beta \wedge \gamma$,
\end{itemize}

and
\begin{itemize}
\item[(ii)] If $\mathfrak{S} \vdash \alpha \Rightarrow \beta$ \, and \, $\mathfrak{S} \vdash \gamma \Rightarrow \delta$ then, \, $\mathfrak{S} \vdash \alpha \wedge \gamma \Rightarrow \beta \wedge \delta$ \, and \,  $\mathfrak{S} \vdash \alpha \vee \gamma \Rightarrow \beta \vee \delta$.
\item[(iii)] $\mathfrak{S} \vdash \alpha \vee (\beta \wedge \gamma) \Leftrightarrow (\alpha \vee \beta) \wedge (\alpha \vee \gamma)$.
\end{itemize}
\end{rem}

\

In what follows, we write $\mathfrak{S} \vdash \Gamma \Leftrightarrow \Sigma$ to indicate that both sequents $\Gamma \Rightarrow \Sigma$ and $\Sigma \Rightarrow \Gamma$ are derivable in $\mathfrak{S}$.

\begin{lem} \label{lem9} Let $\alpha, \beta, \gamma \in Fm$ \, such that \, $\mathfrak{S} \vdash \alpha \Leftrightarrow \beta$. Then,
\begin{itemize}
  \item[(i)] $\mathfrak{S} \vdash \neg \alpha \Leftrightarrow \neg \beta$, 
  \item[(ii)] $\mathfrak{S} \vdash \nabla \alpha \Leftrightarrow \nabla \beta$. 
\end{itemize}
\end{lem}
\begin{dem} (i) is immediate.
(ii):

\begin{prooftree}

\AxiomC{$\alpha \Rightarrow \beta$}

\AxiomC{$\beta \Rightarrow \nabla\beta$}

\LeftLabel{\small(cut)}
\BinaryInfC{$\alpha \Rightarrow \nabla \beta$}
\LeftLabel{\small($\nabla$)}
\UnaryInfC{$\nabla\alpha \Rightarrow \nabla \beta $}
\end{prooftree}

Analogously, it is proved that \, $\nabla\beta \Rightarrow \nabla \alpha$.
\end{dem}

\

Also, since $\neg$ is a De Morgan negation the following result is expected.

\begin{lem} \label{lem8} Let $\alpha,\beta \in Fm$. The following sequents are derivable in $\mathfrak{S}$.
\begin{itemize}
  \item[(i)] $\alpha \Rightarrow \neg \neg \alpha$, \,  $\neg \neg \alpha \Rightarrow \alpha$,
  \item[(ii)] $ \alpha \vee \beta \Rightarrow \neg \neg \alpha \vee \neg \neg\beta$,
  \item[(iii)] $\neg \neg \alpha \wedge \neg \neg\beta \Rightarrow \alpha \wedge \beta$,
  \item[(iv)] $\neg(\alpha \vee \beta) \Rightarrow \neg \alpha \wedge \neg \beta$, 
  \item[(v)] $\neg \alpha \vee \neg \beta \Rightarrow \neg(\alpha \wedge \beta)$,
  \item[(vi)] $\neg \alpha \wedge \neg \beta \Rightarrow \neg(\alpha \vee \beta)$,  
  \item[(vii)] $\neg(\alpha \wedge \beta) \Rightarrow \neg \alpha \vee \neg \beta$,
\end{itemize}
\end{lem}
\begin{dem}(i): Immediate from ($\neg \neg \Rightarrow$) and ($\Rightarrow \neg \neg$).

(ii):

\begin{prooftree}

\AxiomC{$\alpha \Rightarrow \alpha$}
\LeftLabel{\small($\Rightarrow \neg \neg$)}
\UnaryInfC{$\alpha \Rightarrow \neg \neg\alpha$}
\LeftLabel{\small(right w.)}
\UnaryInfC{$\alpha \Rightarrow \neg \neg\alpha, \neg \neg\beta $}
\LeftLabel{\small($\Rightarrow \vee$)}
\UnaryInfC{$\alpha \Rightarrow \neg \neg\alpha\vee \neg \neg\beta$}

\AxiomC{$\beta \Rightarrow \beta$}
\LeftLabel{\small($\Rightarrow \neg \neg$)}
\UnaryInfC{$\beta \Rightarrow \neg \neg\beta$}
\LeftLabel{\small(right w.)}
\UnaryInfC{$\beta \Rightarrow \neg \neg\alpha, \neg \neg\beta $}
\LeftLabel{\small($\Rightarrow \vee$)}
\UnaryInfC{$\beta \Rightarrow \neg \neg\alpha\vee \neg \neg\beta$}
\LeftLabel{\small($\vee\Rightarrow $)}
\BinaryInfC{$\alpha \vee \beta \Rightarrow \neg \neg\alpha\vee \neg \neg\beta$}

\end{prooftree}

\

(iii): Analogous to (ii).

(iv):
\begin{prooftree}
\AxiomC{$\alpha \Rightarrow \alpha$}
\LeftLabel{\small(right w.)}
\UnaryInfC{$\alpha \Rightarrow \alpha,\beta$}
\LeftLabel{\small($\Rightarrow \vee$)}
\UnaryInfC{$\alpha \Rightarrow \alpha \vee \beta $}
\LeftLabel{\small($\neg$)}
\UnaryInfC{$\neg(\alpha \vee \beta)\Rightarrow \neg \alpha $}

\AxiomC{$\beta \Rightarrow \beta$}
\LeftLabel{\small(right w.)}
\UnaryInfC{$\beta \Rightarrow \alpha,\beta$}
\LeftLabel{\small($\Rightarrow \vee$)}
\UnaryInfC{$\beta \Rightarrow \alpha \vee \beta $}
\LeftLabel{\small($\neg$)}
\UnaryInfC{$\neg(\alpha \vee \beta)\Rightarrow \neg \beta$}\LeftLabel{\small($\Rightarrow \wedge$)}
\BinaryInfC{$\neg(\alpha \vee \beta)\Rightarrow \neg \alpha \wedge \neg \beta$}

\end{prooftree}

(v):
\begin{prooftree}
\AxiomC{$\alpha \Rightarrow \alpha$}
\LeftLabel{\small(left w.)}
\UnaryInfC{$\alpha, \beta \Rightarrow \alpha$}
\LeftLabel{\small($\wedge \Rightarrow$)}
\UnaryInfC{$\alpha\wedge \beta \Rightarrow \alpha  $}
\LeftLabel{\small($\neg$)}
\UnaryInfC{$\neg \alpha \Rightarrow \neg (\alpha \wedge \beta) $}

\AxiomC{$\beta \Rightarrow \beta$}
\LeftLabel{\small(left w.)}
\UnaryInfC{$\alpha, \beta \Rightarrow \beta$}
\LeftLabel{\small($\wedge \Rightarrow$)}
\UnaryInfC{$\alpha\wedge \beta \Rightarrow \beta  $}
\LeftLabel{\small($\neg$)}
\UnaryInfC{$\neg \beta \Rightarrow \neg (\alpha \wedge \beta) $}

\LeftLabel{\small($\vee \Rightarrow$)}
\BinaryInfC{$\neg \alpha \vee \neg \beta \Rightarrow \neg (\alpha \wedge \beta)$}
\end{prooftree}

(vi):

\begin{prooftree}
\tiny
\AxiomC{}
\LeftLabel{\tiny(i)}
\UnaryInfC{$\neg \alpha \wedge \neg\beta \Rightarrow \neg  \neg (\neg \alpha \wedge \neg\beta)$}
\LeftLabel{\tiny(r.w.)}
\UnaryInfC{$\neg \alpha \wedge \neg\beta \Rightarrow \neg  \neg (\neg \alpha \wedge \neg\beta), \neg (\alpha\vee \beta)$}

\AxiomC{}
\LeftLabel{\tiny(ii)}
\UnaryInfC{$\alpha \vee \beta \Rightarrow \neg \neg\alpha\vee \neg \neg\beta$}
\LeftLabel{\tiny(r.w.)}
\UnaryInfC{$\alpha \vee \beta \Rightarrow \neg \neg\alpha\vee \neg \neg\beta, \neg(\neg\alpha \vee \neg\beta)$}

\AxiomC{}
\LeftLabel{\tiny(v)}
\UnaryInfC{$\neg \neg\alpha \vee \neg \neg \beta \Rightarrow \neg (\neg\alpha\wedge \neg\beta)$}
\LeftLabel{\tiny(l. w.)}
\UnaryInfC{$\alpha\vee \beta, \neg \neg\alpha \vee \neg \neg \beta \Rightarrow \neg (\neg\alpha\wedge \neg\beta)$} 
\LeftLabel{\tiny(c)}
\BinaryInfC{$\alpha \vee \beta \Rightarrow \neg (\neg \alpha \wedge \neg\beta)$}
\LeftLabel{\tiny($\neg$)}
\UnaryInfC{$\neg \neg (\neg \alpha \wedge \neg\beta) \Rightarrow \neg(\alpha \vee \beta)$}
\LeftLabel{\tiny(l.w.)}
\UnaryInfC{$\neg \alpha \wedge \neg\beta, \neg \neg (\neg \alpha \wedge \neg\beta) \Rightarrow \neg(\alpha \vee \beta)$}

\LeftLabel{\tiny(c)}
\BinaryInfC{$\neg \alpha \wedge \neg\beta \Rightarrow \neg(\alpha \vee \beta)$} 

\end{prooftree}

\

(vii): From (iii), ($\neg$), (i) and the cut rule. 

\end{dem}

\

Finally, the following lemma shows that $\nabla$ enjoys the properties of the respective unary operator of involutive Stone algebras. 

\begin{lem} \label{lem10} Let $\alpha,\beta \in Fm$. Then, 
\begin{itemize}
  \item[(i)] $\mathfrak{S} \vdash \nabla(\alpha \vee \beta) \Leftrightarrow  \nabla \alpha \vee \nabla \beta$,
  \item[(ii)] $\mathfrak{S} \vdash \nabla( \alpha \wedge \beta) \Leftrightarrow \nabla \alpha \wedge \nabla \beta$,
  \item[(iii)] $\mathfrak{S} \vdash \nabla \nabla \alpha  \Leftrightarrow \nabla \alpha$, 
  \item[(iv)] $\mathfrak{S} \vdash \nabla \neg \nabla \alpha  \Leftrightarrow \neg \nabla \alpha$,
	\item[(v)]  $\mathfrak{S} \Rightarrow \neg \nabla \bot$.
\end{itemize}
\end{lem}
\begin{dem}(i)
{\small
\vspace*{-0.5cm}
\begin{prooftree}
\AxiomC{$\alpha \Rightarrow \alpha$}
\UnaryInfC{$\alpha \Rightarrow \alpha, \beta$}
\LeftLabel{\small($\Rightarrow \vee$)}
\UnaryInfC{$\alpha \Rightarrow \alpha \vee \beta$}
\UnaryInfC{$\alpha \Rightarrow \alpha \vee \beta, \nabla(\alpha \vee \beta) $}

\AxiomC{$\alpha \vee \beta\Rightarrow \nabla(\alpha \vee \beta)$}
\UnaryInfC{$\alpha, \alpha \vee \beta\Rightarrow \nabla(\alpha \vee \beta)$}
\LeftLabel{\small(c)}
\BinaryInfC{$\alpha \Rightarrow \nabla(\alpha \vee \beta) $}
\LeftLabel{\small($\nabla$)}
\UnaryInfC{$\nabla \alpha \Rightarrow \nabla(\alpha \vee \beta) $}

\AxiomC{$\beta \Rightarrow \beta$}
\UnaryInfC{$\beta \Rightarrow \alpha, \beta$}
\LeftLabel{\small($\Rightarrow \vee$)}
\UnaryInfC{$\beta \Rightarrow \alpha \vee \beta$}
\UnaryInfC{$\beta \Rightarrow \alpha \vee \beta, \nabla(\alpha \vee \beta) $}

\AxiomC{$\alpha \vee \beta\Rightarrow \nabla(\alpha \vee \beta)$}
\UnaryInfC{$\beta, \alpha \vee \beta\Rightarrow \nabla(\alpha \vee \beta)$}
\LeftLabel{\small(c)}
\BinaryInfC{$\beta \Rightarrow \nabla(\alpha \vee \beta) $}
\LeftLabel{\small($\nabla$)}
\UnaryInfC{$\nabla \beta \Rightarrow \nabla(\alpha \vee \beta) $}

\LeftLabel{\small($\vee \Rightarrow$)}
\BinaryInfC{$\nabla \alpha \vee \nabla \beta \Rightarrow \nabla(\alpha \vee \beta) $}
\end{prooftree}
}

\

\begin{prooftree}
\AxiomC{$\alpha \vee \beta\Rightarrow \alpha \vee \beta$}
\LeftLabel{\small(Lemma \ref{Lem4Equiv})}
\UnaryInfC{$\alpha \vee \beta \Rightarrow \alpha,\beta$}

\AxiomC{$\alpha \Rightarrow \nabla \alpha$}
\LeftLabel{\small(w)}
\UnaryInfC{$\alpha, \beta \Rightarrow \nabla \alpha, \nabla\beta$}

\LeftLabel{\small(cut)}
\BinaryInfC{$\alpha \vee \beta \Rightarrow \nabla \alpha, \nabla\beta$}

\LeftLabel{\small($\nabla$)}
\UnaryInfC{$\nabla(\alpha \vee \beta) \Rightarrow \nabla \alpha, \nabla\beta$}
\LeftLabel{\small($\Rightarrow \vee$)}
\UnaryInfC{$\nabla(\alpha \vee \beta) \Rightarrow \nabla \alpha \vee \nabla\beta$}

\end{prooftree}

\

(ii)

\begin{prooftree}
\AxiomC{$\alpha \Rightarrow \nabla \alpha$}
\LeftLabel{\small(w)}
\UnaryInfC{$\alpha, \beta \Rightarrow  \nabla \alpha$}
\LeftLabel{\small($\wedge \Rightarrow$)}
\UnaryInfC{$\alpha\wedge \beta \Rightarrow  \nabla \alpha$}
\LeftLabel{\small($\nabla$)}
\UnaryInfC{$\nabla(\alpha\wedge \beta)\Rightarrow  \nabla \alpha$}

\AxiomC{$\beta \Rightarrow \nabla \beta$}
\LeftLabel{\small(w)}
\UnaryInfC{$\alpha, \beta \Rightarrow  \nabla \beta$}
\LeftLabel{\small($\wedge \Rightarrow$)}
\UnaryInfC{$\alpha\wedge \beta \Rightarrow  \nabla \beta$}
\LeftLabel{\small($\nabla$)}
\UnaryInfC{$\nabla(\alpha\wedge \beta)\Rightarrow  \nabla \beta$}

\LeftLabel{\small($\Rightarrow \wedge$)}
\BinaryInfC{$\nabla(\alpha\wedge \beta)\Rightarrow \nabla \alpha \wedge \nabla \beta$}

\end{prooftree}

\

\begin{prooftree}
\AxiomC{$\alpha\wedge \beta \Rightarrow \alpha\wedge \beta$}
\LeftLabel{\small(Lemma \ref{Lem4Equiv})}
\UnaryInfC{$\alpha, \beta \Rightarrow  \alpha\wedge \beta$}
\LeftLabel{\small(w)}
\UnaryInfC{$\alpha, \beta \Rightarrow  \alpha\wedge \beta, \nabla(\alpha\wedge \beta)$}

\AxiomC{$\alpha\wedge \beta \Rightarrow \nabla(\alpha\wedge \beta)$}
\LeftLabel{\small(w)}
\UnaryInfC{$\alpha\wedge \beta, \alpha, \beta \Rightarrow \nabla(\alpha\wedge \beta)$}

\LeftLabel{\small(cut)}

\BinaryInfC{$\alpha, \beta \Rightarrow  \nabla(\alpha\wedge \beta)$}
\LeftLabel{\small($\nabla$)}
\UnaryInfC{$\nabla \alpha, \beta \Rightarrow  \nabla(\alpha\wedge \beta)$}
\LeftLabel{\small($\nabla$)}
\UnaryInfC{$\nabla \alpha, \nabla \beta \Rightarrow  \nabla(\alpha\wedge \beta)$}
\LeftLabel{\small($\wedge \Rightarrow$)}
\UnaryInfC{$\nabla \alpha \wedge \nabla \beta \Rightarrow  \nabla(\alpha\wedge \beta)$}
\end{prooftree}

\

(iii)  $\nabla \alpha\Rightarrow \nabla \nabla \alpha$ \, is an instance of the first modal axiom. On the other hand

\

\begin{prooftree}
\AxiomC{$\nabla \alpha \Rightarrow \nabla\alpha$}
\LeftLabel{\small($\nabla$)}
\UnaryInfC{$\nabla \nabla \alpha \Rightarrow \nabla\alpha$}
\end{prooftree}

\

(iv) It is consequence of the first modal axiom and the rule ($\neg \nabla$). 

\

(v)

\begin{prooftree}
\AxiomC{$\bot \Rightarrow $}
\LeftLabel{\small($\nabla$)}
\UnaryInfC{$\nabla \bot \Rightarrow $}
\LeftLabel{\small($\neg$)}
\UnaryInfC{$\Rightarrow \neg \nabla \bot  $}
\end{prooftree}

\end{dem}

\

\section{Inversion principle, soundness and completeness}

In order to prove soundness and completeness of  $\mathfrak{S}$ with respect to {\bf \em Six}, we need to state the notion of {\em validity}.

\begin{defi}\rm
Let $\alpha_1, \dots, \alpha_n, \beta_1, \dots, \beta_m \in Fm$. We shall say that the sequent \, $\alpha_1, \dots, \alpha_n \Rightarrow \beta_1, \dots, \beta_m$ is {\bf valid} iff for all $h\in Hom(Fm,\mathbb{S}_6)$ it holds 

$$h( \bigwedge \limits_{i=1}^{n} \alpha_i) \leq h(\bigvee \limits_{j=1}^{m} \beta_j). $$

If $n = 0$, we take $\bigwedge \limits_{i=1}^{n} \alpha_i = \top$ and if $m=0$, we take $ \bigvee \limits_{j=1}^{m} \beta_j = \bot$.
\end{defi}

From the above definition we have

\begin{prop} $\alpha_1, \dots, \alpha_n \Rightarrow \beta_1, \dots, \beta_m$ is valid \, iff \, $\bigwedge \limits_{i=1}^{n} \alpha_i \models_{\mbox{\bf \em Six }} \bigvee \limits_{j=1}^{m} \beta_j$.
\end{prop}

Besides, we shall say that the rule 

$$\displaystyle\frac{\Gamma_1 \Rightarrow \Sigma_1, \dots, \Gamma_k \Rightarrow \Sigma_k}{\Gamma \Rightarrow \Sigma}$$

{\em preserves validity} if it holds:
 
$$\Gamma_i \Rightarrow \Sigma_i  \mbox{ is valid for } 1\leq i\leq k  \, \,  \mbox{ implies } \, \, \Gamma \Rightarrow \Sigma \mbox{ is valid}. $$

\

Recall that the {\em complexity} of the formula $\alpha$, $c(\alpha)$, is the number of occurrences of connectives in $\alpha$.
We denote $\alpha\mirrormodels\hspace{-.7mm}\models \beta$ to indicate that $\alpha\models_{\mbox{\small \bf \em Six}} \beta$ and $\beta\models_{\mbox{\small \bf \em Six}} \alpha$.
\

\begin{teo}\label{TeoFormaConj} Let $\alpha \in Fm$. Then, \, $\alpha\mirrormodels\hspace{-.7mm}\models\top$ or $\alpha\mirrormodels\hspace{-.7mm}\models\bot$ or

$$\alpha\mirrormodels\hspace{-.7mm}\models \bigwedge \limits_{i=1}^{n} \bigvee \limits_{j=1}^{m_i} \alpha_{ij}$$

where each $\alpha_{ij} \in \{p_k, \neg p_k, \nabla p_k, \nabla \neg p_k, \neg \nabla p_k, \neg \nabla \neg p_k \}$ for some propositional variable $p_k$, $n\geq 0$, $m_i\geq 0$.
\end{teo}
\begin{dem} We use induction on the complexity of $\alpha$.
If $c(\alpha)=0$, then $\alpha$ is $p_k$ or $\top$ or $\bot$. In any of these cases the theorem holds.

(I.H.)Suppose that the theorem holds for any formula with complexity less than $k$, $k\geq 0$. Let $\alpha$ such that $c(\alpha)=k$.

If $\alpha$ is $\alpha^1\wedge \alpha^2$, by  (I.H.), $\alpha^k \mirrormodels\hspace{-.7mm}\models\top$ or $\alpha^k \mirrormodels\hspace{-.7mm}\models\bot$ or

$$\alpha^k \mirrormodels\hspace{-.7mm}\models \bigwedge\limits_{i=1}^{n_k} \bigvee \limits_{j=1}^{m_{k_i}} \alpha^k_{ij}$$ 

for $k=1,2$.

If $\alpha^1 \mirrormodels\hspace{-.7mm}\models\bot$ or $\alpha^2 \mirrormodels\hspace{-.7mm}\models\bot$, then $\alpha \mirrormodels\hspace{-.7mm}\models\bot$.
If $\alpha^1 \mirrormodels\hspace{-.7mm}\models\top$, then $\alpha \mirrormodels\hspace{-.7mm}\models\alpha^2$ and by (I.H.) the theorem holds. Analogously if $\alpha^2 \mirrormodels\hspace{-.7mm}\models\top$. 

Otherwise, $\alpha^1 \mirrormodels\hspace{-.7mm}\models  \bigwedge\limits_{i=1}^{n_1} \bigvee \limits_{j=1}^{m_{1_i}} \alpha^1_{ij}$ and $\alpha^2 \mirrormodels\hspace{-.7mm}\models  \bigwedge\limits_{i=1}^{n_2} \bigvee \limits_{j=1}^{m_{2_i}} \alpha^2_{ij}$. Then, $\alpha=\alpha^1\wedge \alpha^2 \mirrormodels\hspace{-.7mm}\models \bigwedge\limits_{i=1}^{n_1} \bigvee \limits_{j=1}^{m_{1_i}} \alpha^1_{ij} \wedge \bigwedge\limits_{i=1}^{n_2} \bigvee \limits_{j=1}^{m_{2_i}} \alpha^2_{ij}$ and this last expression is a conjunction of terms which are disjunction of formulas whose forms are in  $\{p_k, \neg p_k, \nabla p_k, \nabla \neg p_k, \neg \nabla p_k, \neg \nabla \neg p_k \}$.\\[2mm]
If $\alpha$ is $\alpha^1\vee \alpha^2$ with $\alpha^1$ and $\alpha^2$ as in the previous case. Then, if $\alpha^1 \mirrormodels\hspace{-.7mm}\models\top$ or $\alpha^2 \mirrormodels\hspace{-.7mm}\models\top$, then $\alpha \mirrormodels\hspace{-.7mm}\models\top$.
If $\alpha^1 \mirrormodels\hspace{-.7mm}\models\bot$, then $\alpha \mirrormodels\hspace{-.7mm}\models\alpha^2$ and by (I.H.) the theorem holds. Analogously if $\alpha^2 \mirrormodels\hspace{-.7mm}\models\bot$.
Otherwise,  $\alpha=\alpha^1\vee \alpha^2 \mirrormodels\hspace{-.7mm}\models \bigwedge\limits_{i=1}^{n_1} \bigvee \limits_{j=1}^{m_{1_i}} \alpha^1_{ij} \vee \bigwedge\limits_{i=1}^{n_2} \bigvee \limits_{j=1}^{m_{2_i}} \alpha^2_{ij}$. By Corollary \ref{CorEqSemIgualdad} and having in mind that every involutive Stone algebra is, in particular, a distributive lattice,  using distributivity as many times as necessary we can transform the expression  $\bigwedge\limits_{i=1}^{n_1} \bigvee \limits_{j=1}^{m_{1_i}} \alpha^1_{ij} \vee \bigwedge\limits_{i=1}^{n_2} \bigvee \limits_{j=1}^{m_{2_i}} \alpha^2_{ij}$ into one expression  of the form \, $\bigwedge \bigvee \beta_t$ \, where \, $\beta_t\in \{p_k, \neg p_k, \nabla p_k, \nabla \neg p_k, \neg \nabla p_k, \neg \nabla \neg p_k \}$ for some propositional variable $p_k$.\\[2mm]
If $\alpha$ is $\neg \alpha^1$, by (I.H.), $\alpha^1 \mirrormodels\hspace{-.7mm}\models\top$ or $\alpha^1 \mirrormodels\hspace{-.7mm}\models\bot$ or $\alpha^1 \mirrormodels\hspace{-.7mm}\models \bigwedge\limits_{i=1}^{n_1} \bigvee \limits_{j=1}^{m_{1_i}} \alpha^1_{ij}$. If $\alpha^1 \mirrormodels\hspace{-.7mm}\models\top$, then $\alpha \mirrormodels\hspace{-.7mm}\models\bot$. If $\alpha^1 \mirrormodels\hspace{-.7mm}\models\bot$, then $\alpha \mirrormodels\hspace{-.7mm}\models\top$. Otherwise, by Lemma \ref{lem9} and Remark \ref{Obs1}, $\alpha=\neg \alpha^1 \mirrormodels\hspace{-.7mm}\models \neg \bigwedge\limits_{i=1}^{n_1} \bigvee \limits_{j=1}^{m_{1_i}} \alpha^1_{ij}$. Since every involutive Stone algebra is, in particular, a De Morgan algebra we have $\alpha\mirrormodels\hspace{-.7mm}\models \neg \bigwedge\limits_{i=1}^{n_1} \bigvee \limits_{j=1}^{m_{1_i}} \alpha^1_{ij}\mirrormodels\hspace{-.7mm}\models \bigvee\limits_{i=1}^{n_1} \bigwedge \limits_{j=1}^{m_{1_i}} \neg \alpha^1_{ij}$.
Again, by distributivity and (DM1), $\neg \neg x \approx x$, we have that the expression  $\bigvee\limits_{i=1}^{n_1} \bigwedge \limits_{j=1}^{m_{1_i}} \neg \alpha^1_{ij}$ is equivalent to one of the type  \, $\bigwedge \bigvee \beta_t$ \, where each \, $\beta_t\in \{p_k, \neg p_k, \nabla p_k, \nabla \neg p_k, \neg \nabla p_k, \neg \nabla \neg p_k \}$ for some propositional variable $p_k$.\\[2mm]

Finally, if $\alpha$ is $\nabla \alpha^1$ and $\alpha^1 \mirrormodels\hspace{-.7mm}\models\top$, then $\alpha \mirrormodels\hspace{-.7mm}\models\top$. If $\alpha^1 \mirrormodels\hspace{-.7mm}\models\bot$, then $\alpha \mirrormodels\hspace{-.7mm}\models\bot$. Otherwise, $\alpha^1 \mirrormodels\hspace{-.7mm}\models \bigwedge\limits_{i=1}^{n_1} \bigvee \limits_{j=1}^{m_{1_i}} \alpha^1_{ij}$ and from Lemma \ref{propcuadrado} (IS7), (IS8) and $\nabla$ being a lattice--homomorphism we have $\alpha=\nabla \alpha^1 \mirrormodels\hspace{-.7mm}\models  \bigwedge\limits_{i=1}^{n_1} \bigvee \limits_{j=1}^{m_{1_i}} \nabla\alpha^1_{ij}$,
where $\nabla\alpha^1_{ij}$ has the form of some element of  \\
$\{p_k, \neg p_k, \nabla p_k, \nabla \neg p_k, \neg \nabla p_k, \neg \nabla \neg p_k \}$ for some propositional variable $p_k$.
\end{dem}

\

\begin{rem} From Theorem \ref{TeoFormaConj} we have that every formula  $\alpha \in Fm$ is logically equivalent (in {\bf Six}) to  $\top$  or $\bot$ or to a formula of the form  $\bigwedge \limits_{i=1}^{n} \bigvee \limits_{j=1}^{m_i} \alpha_{ij}$. If $\alpha$ is not equivalent to $\bot$ nor to  $\top$, we shall say that $\bigwedge \limits_{i=1}^{n} \bigvee \limits_{j=1}^{m_i} \alpha_{ij}$ is one {\bf conjunctive form} of $\alpha$. Every formula has infinite conjunctive forms. We shall denote by $\bar{\alpha}$ any conjunctive form of $\alpha$. If $\alpha$ is equivalent to $\bot$ or  $\top$, then  $\bot$ ( $\top$) is a conjunctive form of $\alpha$. Besides, we say that the $\alpha_{ij}$'s are the {\em blocks} of  $\bar{\alpha}$.
\end{rem}

If $\alpha$ is a conjunctive form, we denote by $\# \alpha$ the number of blocks of $\alpha$. That is, if $\alpha$ is $\bot$ or  $\top$, then $\# \alpha = 0$. Otherwise, $\alpha$ is  $\bigwedge \limits_{i=1}^{n} \bigvee \limits_{j=1}^{m_i} \alpha_{ij}$ and then $\# \alpha = m_1+\dots +m_n$.

\

\begin{eje} Consider the formula \, $\nabla((p_1\wedge \neg \nabla p_2) \vee \nabla p_2)$. By Lemma \ref{lem10} (i), we have

$$\nabla((p_1\wedge \neg \nabla p_2) \vee \nabla p_2) \, \Leftrightarrow \, \nabla(p_1\wedge \neg \nabla p_2) \vee \nabla \nabla p_2$$

and, by Lemma \ref{lem10} (ii) and  Remmark \ref{Obs1} (i),

$$\nabla(p_1\wedge \neg \nabla p_2) \vee \nabla \nabla p_2 \, \Leftrightarrow \, (\nabla p_1\wedge \nabla \neg \nabla p_2) \vee \nabla \nabla p_2$$

Now, $\nabla p_1\wedge \nabla \neg \nabla p_2 \, \Leftrightarrow \, \nabla p_1\wedge \neg \nabla p_2$, by Remark \ref{Obs1} (i) and Lemma \ref{lem10} (iv). Then, by Remark \ref{Obs1} (i)

$$(\nabla p_1\wedge \nabla \neg \nabla p_2) \vee \nabla \nabla p_2 \, \Leftrightarrow \, (\nabla p_1 \wedge \neg \nabla p_2) \vee \nabla \nabla p_2$$

From Lemma \ref{lem10} (iv) and  Remark \ref{Obs1} (i), we have

$$(\nabla p_1\wedge  \neg \nabla p_2) \vee \nabla \nabla p_2 \, \Leftrightarrow \, (\nabla p_1 \wedge \neg \nabla p_2) \vee \nabla p_2$$

Finally, from Remark \ref{Obs1} (iii), we have

$$(\nabla p_1 \wedge \neg \nabla p_2) \vee \nabla p_2  \, \Leftrightarrow \, (\nabla p_1  \vee \nabla p_2) \wedge (\neg \nabla p_2 \vee \nabla p_2) $$

This last formula is a conjunctive form logically equivalent to the original one, and we have proved that $\nabla((p_1\wedge \neg \nabla p_2) \vee \nabla p_2)$ and $(\nabla p_1  \vee \nabla p_2) \wedge (\neg \nabla p_2 \vee \nabla p_2)$ are inter--derivable. That is, the following sequents are derivable in $\mathfrak{S}$

$$\nabla((p_1\wedge \neg \nabla p_2) \vee \nabla p_2) \, \Leftrightarrow \, (\nabla p_1  \vee \nabla p_2) \wedge (\neg \nabla p_2 \vee \nabla p_2)$$

\

Besides, $\# (\nabla p_1  \vee \nabla p_2) \wedge (\neg \nabla p_2 \vee \nabla p_2) = 4$
\end{eje}

\

\begin{lem}  \label{lemFormaConj} Let $\alpha\in Fm$. Then, there exists a conjunctive form $\bar{\alpha}$ of  $\alpha$ such that \, $\alpha \Leftrightarrow \bar{\alpha}$.
\end{lem}
\begin{dem} The proof is by induction on the complexity of $\alpha$ and using Remark \ref{Obs1} and Lemmas \ref{lem9}, \ref{lem8}, \ref{lem10}. We leave the details to the patient reader.
\end{dem}

 
\

\begin{teo}\label{teoCorrectitud}(Soundness) Every sequent $\Gamma \Rightarrow \Sigma$ derivable in $\mathfrak{S}$ is valid.
\end{teo}
\begin{dem} It is easy to check that the axioms are valid sequents and that the structural rules preserve validity.
On the other hand, taking into account that $\mathbb{S}_{6}$ is a De Morgan algebra, (DM1), (IS7) and (IS8) we see clearly that  ($\neg$), ($\neg \neg \Rightarrow$), ($\Rightarrow \neg \neg$) and $(\neg\nabla)$ preserve validity.

Finally, let's see that 

$$ \mbox{($\nabla$) \, } \displaystyle \frac{\Gamma, \alpha \Rightarrow \nabla \Sigma } {\Gamma, \nabla \alpha \Rightarrow \nabla \Sigma} $$

preserves validity. Suppose that (1) $\bigwedge \limits_{i=1}^{n} h(\alpha_i) \wedge h(\alpha) \leq \bigvee \limits_{j=1}^{m} h(\nabla\beta_j)$. Then, since $\nabla$ is a lattice--homomorphism, \, $\bigwedge \limits_{i=1}^{n} h(\alpha_i) \wedge h(\alpha) \leq \bigvee \limits_{j=1}^{m} \nabla h(\beta_j) = \nabla \left( \bigvee \limits_{j=1}^{m} h(\beta_j) \right)$. However, in $\mathbb{S}_6$ we have that $\nabla x\in \{0,1\}$ for all $x$. Then, $\nabla \left( \bigvee \limits_{j=1}^{m} h(\beta_j) \right)\in \{0,1\}$. If 
$\nabla \left( \bigvee \limits_{j=1}^{m} h(\beta_j) \right)=0$, then, by (1),  $\bigwedge \limits_{i=1}^{n} h(\alpha_i) \wedge h(\alpha)=0$ and then 
$\bigwedge \limits_{i=1}^{n} h(\alpha_i)=0$ or $ h(\alpha)=0$. In both cases we have $\bigwedge \limits_{i=1}^{n} h(\alpha_i) \wedge h(\nabla \alpha)=0$ and therefore $\bigwedge \limits_{i=1}^{n} h(\alpha_i) \wedge h(\nabla \alpha)\leq \bigvee \limits_{j=1}^{m}  h(\nabla \beta_j)$. If $\nabla \left( \bigvee \limits_{j=1}^{m} h(\beta_j) \right)=1$, trivially we have $\bigwedge \limits_{i=1}^{n} h(\alpha_i) \wedge h(\nabla \alpha)\leq \bigvee \limits_{j=1}^{m}  h(\nabla \beta_j)$.
\end{dem}

\

Next, we shall see that in $\mathfrak{S}$ the {\em Inversion principle} holds.

\

\begin{lem} \label{lemInv} (Inversion principle) For every rule \, $\displaystyle\frac{S_1,\dots,S_k}{S}$ \, of $\mathfrak{S}$, except left and right weakening, it holds:
if $S$ is valid, then $S_i$ is valid for $i\in\{1, \dots, k\}$.
\end{lem}
\begin{dem} Since rules ($\wedge \Rightarrow$), ($\Rightarrow\wedge$), ($\vee \Rightarrow$) and ($\Rightarrow \vee$) are the same as in the classical case; and taking into account that involutive Stone algebras are distributive lattices, the principle holds for them as in the classical case.
Since every involutive Stone algebra is a De Morgan algebra, rules ($\neg$), ($\neg \neg \Rightarrow$) and ($\Rightarrow \neg \neg$) verify the principle.
Rule ($\neg\nabla\Rightarrow$) verifies the principle by (IS8).

\

Let's see that ($\nabla$) verifies the principle. Suppose that the sequent \, \,  $\Gamma, \nabla \alpha \Rightarrow \nabla \Sigma$ is valid. Then, for every $h\in Hom(Fm,\mathbb{S}_6)$ we have that 

$$h(\bigwedge \limits_{\gamma \in \Gamma}\gamma) \wedge \nabla h(\alpha) = h( \bigwedge \limits_{\gamma \in \Gamma}\gamma \wedge \nabla \alpha) \leq h(\bigvee \limits_{\beta \in \Sigma} \nabla \beta) =  h(\nabla \bigvee \limits_{\beta \in \Sigma} \beta) = \nabla h(\bigvee \limits_{\beta \in \Sigma} \beta)$$

By the definition of $\nabla$ in $\mathbb{S}_6$ we have that $\nabla h(\bigvee \limits_{\beta \in \Sigma} \beta) \in \{0,1\}$. If $\nabla h(\bigvee \limits_{\beta \in \Sigma} \beta)=0$, then $h(\bigwedge \limits_{\gamma \in \Gamma}\gamma) \wedge \nabla h(\alpha)=0$. In $\mathbb{S}_6$, we can assert that $h(\bigwedge \limits_{\gamma \in \Gamma}\gamma) =0$ or  $\nabla h(\alpha)=0$. If $h(\bigwedge \limits_{\gamma \in \Gamma}\gamma) =0$ then $h( \bigwedge \limits_{\gamma \in \Gamma}\gamma \wedge \alpha) \leq h(\bigvee \limits_{\beta \in \Sigma} \nabla \beta)$. If $\nabla h(\alpha)=0$, then $h(\alpha)=0$ and therefore $h(\bigwedge \limits_{\gamma \in \Gamma}\gamma) \wedge h(\alpha)=0$ and $h( \bigwedge \limits_{\gamma \in \Gamma}\gamma \wedge \alpha) \leq h(\bigvee \limits_{\beta \in \Sigma} \nabla \beta)$.

On the other side, if $\nabla h(\bigvee \limits_{\beta \in \Sigma} \beta)=1$, clearly  $h( \bigwedge \limits_{\gamma \in \Gamma}\gamma \wedge \alpha) \leq h(\bigvee \limits_{\beta \in \Sigma} \nabla \beta)$.
\end{dem}

\

\begin{teo}\label{teoCompletitud} (Completeness) Every sequent valid in $\mathfrak{S}$ is derivable in $\mathfrak{S}$.
\end{teo}
\begin{dem} We show that every valid sequent $\Gamma \Rightarrow \Sigma$ has a $\mathfrak{S}$ proof, by induction on the total number of logical
connectives $\vee$ and $\wedge$ occurring in $\Gamma \Rightarrow \Sigma$. In virtue of Lemma \ref{lemFormaConj}, we assume that every formula in $\Gamma \cup \Sigma$ is in conjunctive form. 

For the base case, every formula
in $\Gamma$ and $\Sigma$ is a block in $\{p, \neg p, \nabla p, \nabla \neg p, \neg \nabla p, \neg \nabla \neg p \}$ for some propositional variable $p$ or one of the constants $\bot$ or $\top$. Since the sequent is valid, we have the following cases:

\begin{itemize}
  \item[](1) $\bot \in \Gamma$ or $\top \in \Sigma$,  
  \item[](2.1) $p \in \Gamma$ and $p \in \Sigma$; (2.2) $p \in \Gamma$ and  $\nabla p \in \Sigma$,
  \item[](3.1) $\neg p \in \Gamma$ and $\neg p \in \Sigma$; (3.2) $\neg p \in \Gamma$ and  $\nabla \neg p \in \Sigma$,
  \item[](4.1) $\nabla p \in \Gamma$ and $\nabla p \in \Sigma$; (4.2) $\nabla \neg p \in \Gamma$ and  $\nabla \neg p \in \Sigma$,
  \item[](5.1) $\neg \nabla \neg p \in \Gamma$ and $\neg \nabla \neg p \in \Sigma$; (5.2) $\neg \nabla \neg p \in \Gamma$ and  $p \in \Sigma$,
  \item[](6) $\nabla p, \neg \nabla  p \in \Sigma$.
\end{itemize}

In case (1), $\Gamma \Rightarrow \Sigma$ can be derived from ($\bot$) or ($\top$) and by weakenings. In cases (2.1), (3.1), (4.1), (4.2) and (5.1), $\Gamma \Rightarrow \Sigma$ can be derived from the structural axiom and by weakenings. In cases (2.2), (3.2),  $\Gamma \Rightarrow \Sigma$, can be derived from the first modal axiom and by weakenings. In case (5.2), we can derive $\neg \nabla \neg p \Rightarrow p$ as follows

\begin{prooftree}
\AxiomC{$\neg p \Rightarrow \nabla \neg p$}
\LeftLabel{\small($\neg$)}
\UnaryInfC{$\neg \nabla \neg p \Rightarrow \neg \neg p$}

\AxiomC{$\neg \neg p \Rightarrow  p$}
\LeftLabel{\small(cut)}
\BinaryInfC{$\neg \nabla \neg p \Rightarrow p$}
\end{prooftree}

and from this we derive $\Gamma \Rightarrow \Sigma$ by weakenings. In case (6), we derive $\Gamma \Rightarrow \Sigma$ by using the second modal axiom, Lemma \ref{Lem4Equiv} and weakenings.

For the induction step, let $\alpha$ be any formula which is not a block and
not a constant in $\Gamma$ or $\Sigma$. Then by the definition of propositional formula
$\alpha$ must have one of the forms ($\beta  \vee \gamma$) or ($\beta  \wedge \gamma$). If $\alpha$ has the form ($\beta  \vee \gamma$) and $\nabla p$ is a block of $\beta$ and $\neg \nabla p$ is a block of $\gamma$, then $\Gamma \Rightarrow \Sigma$ by using the second modal axiom, Lemma \ref{Lem4Equiv} and weakenings.

Otherwise $\Gamma \Rightarrow \Sigma$ can be derived from $\vee$-- introduction or $\wedge$--introduction, respectively, using either the left case or the right case, depending on whether $\alpha$ is in $\Gamma$ or $\Sigma$, but not using weakenings. In each case,
each top sequent of the rule will have at least one fewer connective $\wedge$ or $\vee$ than
$\Gamma \Rightarrow \Sigma$, and the sequent is valid by the inversion principle. Hence each
top sequent has a derivation in $\mathfrak{S}$, by the induction hypothesis. 

\end{dem}

\section{Conclusions}

This paper was dedicated to studying the logic that preserves degrees of truth with respect to involutive Stone algebras, named {\bf \em Six}. We showed that this is a six-valued logic that can be defined in terms of four matrices.  {\bf \em Six} turns out to be to be an interesting logic which combines nice characteristics of Belnap's four-valued logic $FOUR$ and $n$-valued \L ukasiewicz logics ($2\leq n\leq5$). As a paraconsistent logic, it is a genuine LFI with a consistency operator that distinguishes between consistency and contradiction. Finally, we provide a syntactic presentation of {\bf \em Six} by means of the sequent calculus $\mathfrak{S}$. We think that $\mathfrak{S}$ deserves a deeper study, for instance, it is an open question weather it admits cut elimination. We presume this is not the case in view of the unusual form of the logical rules that regulate the behavior of the modal operator $\nabla$. In this sense, we think that the framework of hypersequents would be a more suitable tool for presenting a Gentzen--style system for {\bf \em Six} with the cut--elimination property.

\noindent {\bf Acknowledgements:}  We would like to thank the anonymous referees for their helpful comments and suggestions which allowed us to improve the final version of our paper.

\end{document}